\documentclass{article}

\usepackage[round, authoryear]{natbib}
\usepackage{graphicx}
\usepackage{amsmath,amsthm,amsfonts,amssymb,eucal,eufrak}
\usepackage{color}

\usepackage{fullpage}
\usepackage{setspace}

\usepackage{hyperref} 
\hypersetup{backref, colorlinks=true, citecolor=blue, linkcolor=blue}
\newcommand{\rref}[1]{\hyperref[#1]{\ref*{#1}}}

\newtheorem{thm}{Theorem}[section]

\newtheorem{cor}[thm]{Corollary}
\newtheorem{lem}[thm]{Lemma}

\def\OO#1{\mbox{\rm O}\left(#1\right)}

\title{Autoregressive Process Modeling via the Lasso Procedure}
\author{Yuval Nardi\thanks{Email: {\tt yuval@stat.cmu.edu}}\\
Department of Statistics\\
Carnegie Mellon University\\
Pittsburgh, PA 15213-3890 USA
\and
Alessandro Rinaldo\thanks{Email: {\tt arinaldo@stat.cmu.edu}}\\
Department of Statistics\\
Carnegie Mellon University\\
Pittsburgh, PA 15213-3890 USA}

\date{}

\begin{document}

 \nocite{*}
\maketitle

\begin{abstract}
The Lasso is a popular model selection and estimation procedure for linear models that enjoys nice theoretical properties. In this paper, we study the Lasso estimator for fitting autoregressive time series models. We adopt a double asymptotic framework where the maximal lag may increase with the sample size. We derive theoretical results establishing various types of consistency. In particular, we derive conditions under which  %that, under suitable conditions, and with an appropriate choice of the autoregressive order, 
the Lasso estimator for the autoregressive coefficients is model selection consistent, estimation consistent and prediction consistent. Simulation study results are reported.
\end{abstract}

\section{Introduction}\label{sec:intro}

Classical stationary time series modeling assumes that data are a realization of a mix of autoregressive processes and moving average processes, or an ARMA model \citep[see, e.g.][]{MR1093459}. Typically, both estimation and model fitting rely on the assumption of fixed and low dimensional parameters and include  $(i)$ the estimation of the appropriate coefficients under the somewhat unrealistic assumption that the orders of the AR and of the MA processes are known in advance, or $(ii)$ some model selection procedures that sequentially fit models of  increasing dimensions. In practice, however, it is very difficult to verify the assumption that the realized series does come from an ARMA process. Instead, it is usually assumed that the given data are a realization of a \emph{linear} time series, which may be represented by an infinite-order autoregressive process. Some study has been done on the accuracy of an AR approximation for these processes: see \cite{Shibata}, \cite{zeevi} and \cite{Ing}. In particular, \citet{zeevi} propose a nonparametric minimax approach and assess the accuracy of a finite order AR process in terms of both estimation and prediction. 

This paper is concerned with fitting autoregressive time series models with the Lasso.  The Lasso procedure, proposed originally by \citet{lasso}, is one of the most popular approach for model selection in linear and generalized linear models, and has been studied in much of the recent literature; see, e.g., \cite{fan-2004-32}, \cite{zhaobin}, \cite{adaptive}, \cite{wainwright}, \cite{LaffertySpam}, and \cite{NardiRinaldo}, to mention just a few.
The Lasso procedure has the advantage of simultaneously performing model selection and estimation, and has been shown to be effective even in high dimensional settings where the dimension of the parameter space grows with the sample size $n$. In the context of an autoregressive modeling, the Lasso features become especially advantageous, as both the AR order, and the corresponding AR coefficients can be estimated simultaneously. %The Lasso procedure is one of the most popular approach for model selection, and is the subject of many current studies; see, e.g., \cite{fan-2004-32}, \cite{zhaobin}, \cite{adaptive}, \cite{wainwright}, \cite{LaffertySpam}, and \cite{NardiRinaldo}, to mention just a few. 
\citet{MR2301500} study linear regression with autoregressive errors. They adapt the Lasso procedure to shrink both the regression coefficients and the autoregressive coefficients, under the assumption  that the autoregressive order is fixed.

For the autoregressive models we consider in this work, the number of parameters, or equivalently, the maximal possible lag, grows with the sample size. We refer to this scheme as a double asymptotic framework. The double asymptotic framework enables us to treat the autoregressive order as virtually infinite. The autoregressive time series with an increasing number of parameters lies between a fixed order AR time series and an infinite-order AR time series. This limiting process belongs to a family which is known to contain many ARMA processes \citep[see][]{zeevi}. In this paper we show that the  Lasso procedure is particularly adequate for this double asymptotic scheme. %and is known to produce good results in many settings (see the above mentioned papers).  

The rest of the paper is organized as follows. The next section formulates the autoregressive modeling scheme and defines the Lasso estimator associated with it. Asymptotic properties of the Lasso estimator are presented in Section \ref{sec:results}. These include model selection consistency (Theorem \ref{thm:sign}), estimation consistency (Theorem \ref{thm:estimation}), and prediction consistency (Corollary \ref{cor:prediction}). Proofs are deferred to Section \ref{sec:proofs}. A simulation study, given in Section \ref{sec:simul}, accompany the theoretical results. Discussion and concluding remarks appear in Section \ref{sec:discussion}.

\section{Penalized autoregressive modeling}\label{sec:ar}

In this section we describe our settings and set up the notation.

%In this section we set-up the notation and describe the AR model along with the Lasso estimator for the AR coefficients to be considered.

We assume that $X_1, \ldots, X_n$ are $n$ observations from an AR($p$) process:
\begin{equation}\label{eq:ar_model}
X_t = \phi_1 X_{t-1}+\ldots +\phi_p X_{t-p}+Z_t \quad ,\quad t=1,\ldots, n \; ,
\end{equation}
where $\{Z_t\}$ is a random sequence of independent Gaussian variables  with $\mathbb{E}Z_t=0$, $\mathbb{E}|Z_t|^2=\sigma^2$ and $\text{cov}(Z_t,X_s)=0$ for all $s<t$. The last requirement is standard, and rely on a reasoning under which the process $\{X_t\}$  does not depend on future values of the driving Gaussian noise. The assumption about Gaussianity of $\{Z_t\}$ is by no means necessary, and can be relaxed. It does, however, facilitate our theoretical investigation and the presentation of various results, and therefore, it is in effect throughout the article.  In Section \ref{sec:discussion} we comment on how to modify our assumptions and proofs to allow for non-Gaussian innovations $\{Z_t\}$. 

We further assume that $\{ X_t\}$ is \emph{causal}, meaning that there exists a sequence of constants $\{\psi_j\}$, $j=0,1, \ldots$, with absolutely convergent series, $\sum_{j=0}^\infty |\psi_j|<\infty$, such that $\{X_t\}$ has a MA($\infty$) representation: 
\begin{equation}\label{eq:ma_infty}
X_t = \sum_{j=0}^\infty \psi_jZ_{t-j} \; ,
\end{equation}
the series being absolutely convergent with probability one. Equivalently, we could stipulate that $\{X_t\}$ is purely non-deterministic, and then obtain representation (\ref{eq:ma_infty}), with $\psi_0=1$ and $\sum_{j=0}^\infty \psi_j^2<\infty$, directly from the Wold decomposition \citep[see, e.g.][]{MR1093459}. A necessary and sufficient condition for causality is that $1-\phi_1 z-\ldots -\phi_p z^p \ne 0$ for all complex $z$ within the unit disc, $|z|\leq 1$. Notice that causality of $\{X_t\}$, and Gaussianity of $\{Z_t\}$, together imply Gaussianity of $\{X_t\}$. This follows from the fact that mean square limits of Gaussian random variables are again Gaussian. The mean and variance of $X_t$ are given, respectively, by $\mathbb{E}X_t=0$, $\mathbb{E}|X_t|^2=\sigma^2\sum_{j=0}^\infty \psi_j^2$. We assume, for simplicity, and without any loss of generality, that $\mathbb{E}|X_t|^2=1$, so that $\sum_{j=0}^\infty \psi_j^2=\sigma^{-2}$. Let $\gamma(\cdot)$ be the autocovariance function given by  $\gamma(k)=\mathbb{E}X_tX_{t+k}$, and let $\Gamma_p=\big(\gamma(i-j)\big)_{i,j=1,\ldots, p}$, the $p\times p$ autocovariance matrix, of lags smaller or equal to $p-1$.

We now describe the penalized $\ell_1$ least squares estimator of the AR coefficients. Let $y=(X_1, \ldots, X_n)'$, $\phi=(\phi_1, \ldots, \phi_p)'$, and $Z=(Z_1, \ldots, Z_n)'$, where apostrophe denotes transpose. Define the $n\times p$ matrix $X$ with entry $X_{t-j}$ in the $t$th row and $j$th column, for $t=1,\ldots, n$ and $j=1,\ldots, p$. The Lasso-type estimator $\hat\phi_n \equiv \hat \phi_n(\Lambda_n)$ is defined to be the minimizer of:
\begin{equation}\label{eq:lasso_ar}
\frac{1}{2n}\|y-X\phi\|^2
+\lambda_n \sum_{j=1}^p \lambda_{n,j}|\phi_j|  \; ,
\end{equation}
where $\Lambda_n=\{\lambda_n, \{\lambda_{n,j}\, ,\, j=1,\ldots, p\}\}$ are tuning parameters, and $\|\cdot\|$ denotes the $l_2$-norm. Here, $\lambda_n$ is a grand tuning parameter, while the $\{\lambda_{n,j}, j =1, \ldots, p\}$ are specific tuning parameters associated with predictors $X_{t-j}$. 
%In our objective function  (\ref{eq:lasso_ar}), 
The Lasso solution (\ref{eq:lasso_ar}) will be sparse, as some of the autoregressive coefficients will be to set to (exactly) zero, depending on the choice tuning parameters $\Lambda_n$.  Naturally, one may want to further impose that $\lambda_{n,j}<\lambda_{n,k}$ for lags values satisfying $j<k$, to encourage even sparser solutions, although this is not assumed throughout. 
The idea of using $\ell_1$ regularization to penalize differently the model parameters, as we do in (\ref{eq:lasso_ar}), was originally proposed by \citet{adaptive} under the name of adaptive Lasso. As shown in \citet{adaptive}, from an algorithmic  point of view, the solution to our adaptive Lasso (\ref{eq:lasso_ar})  can be obtained by a slightly modified version of the LARS algorithm  of \citet{lars}. A possible choice for $\lambda_{n,j}$ would be to use the inverse least squares estimates, as in \citet{adaptive}, but this is not pursued here.

%The idea of using  different tuning parameters for every covariate, which w, appears in other contexts.  \citet{adaptive} uses similar set of tuning parameters, but picks the $\lambda_{n,j}$s adaptively using the ordinary least squares estimator. It is shown there that the weighted (adaptive) Lasso estimator can be solved by means of the LARS algorithm (\cite{lars}). This is, clearly, true also for the present problem. \citet{yuan} define the group lasso estimator and use the dimension of sub-groups as weights $\lambda_{n,j}$.

As mentioned before, we consider a double asymptotic framework, in which the number of parameters $p\equiv p_n$ grows with $n$ at a certain rate. Clearly, the ``large $p$ small $n$'' ($p\gg n$) scenario, which is an important subject of many of nowadays articles, is not adequate here. Indeed, one might be suspicious about the statistical properties of the proposed estimator even when $p$ is comparable with $n$ ($p<n$, but is close to $n$). Accounting for the mechanism of the autoregressive progress, one is led to think that $p$ should grow with $n$ at a much slower rate. This article shows that the choice of $p=\OO{\log n}$ will lead to nice asymptotic properties of the proposed procedure (\ref{eq:lasso_ar}). Such a choice of the AR order arises also in \citet{zeevi}, who prove minimax optimality for a different regularized least squares estimator. Moreover, as pointed out in \cite{zeevi}, the same order of $p$ arises also in spectral density estimation (see \cite{efromovich}). Finally, similar rate appears also, in a different context, in \cite{bickel_levina}.

In classical linear time series modeling, one usually  attempts to  fit sequentially an AR($p$) with increasing orders of the maximal lag $p$ (or by fixing $p$ and then estimating the coefficients). The Lasso-type estimator of scheme (\ref{eq:lasso_ar}) will shrink down to zero irrelevant predictors. Thus, not only that model selection and estimation will occur simultaneously, but the fitted (selected) model will be chosen among all relevant AR($p$) processes, with $p=\OO{\log n}$.

\section{Asymptotic Properties of the Lasso}\label{sec:results}

In this section we derive the asymptotic properties of the Lasso estimator $\hat \phi_n$. These include model selection consistency, estimation consistency and prediction consistency. We briefly describe each type of consistency, develop the needed notation, and present the results, with proofs relegated to Section \ref{sec:proofs}.

\subsection{Model Selection Consistency}

We assume that the AR($p$) process (\ref{eq:ar_model}) is generated according to a true, unknown parameter $\phi^*=(\phi^*_1,\ldots, \phi^*_p)$. When $p$ is large, it is not unreasonable to believe that this vector is sparse, meaning that only a subset of potential predictors are relevant. Model selection consistency is about recovering the sparsity structure of the true, underlying parameter $\phi^*$.

For any vector $\phi\in\mathbb{R}^p$, let $\text{sgn}(\phi)=(\text{sgn}(\phi_1), \ldots, \text{sgn}(\phi_p))$, where $\text{sgn}(\phi_j)$ is the sign function taking values $-1, 0$ or $ 1$, according to as $\phi_j<0, \phi_j=0$ or $\phi_j>0$, respectively. A given estimator $\hat\phi_n$ is said to be \emph{sign consistent} if $\text{sgn}(\hat\phi_n)=\text{sgn}(\phi^*)$, with probability tending to one, as $n$ tends to infinity, i.e.,
\begin{equation}\label{eq:sign_consistency}
\mathbb{P}(\text{sgn}(\hat\phi_n) =\text{sgn}(\phi^*))\longrightarrow 1 \qquad ,\qquad  n\rightarrow\infty \; .
\end{equation}
Let $S=\{j  \, : \, \phi_j^*\ne 0\}=\text{supp}(\phi^*)\subset\{1,2,\ldots, p\}$. A weaker form of model selection consistency, implied by the sign consistency, only requires that, with probability tending to $1$, $\phi^*$ and $\hat\phi_n$ have the same support. %As noted in \cite{zhaobin}, sign consistency excludes situations for which the signs of the estimated are matched incorrectly, as it may be the case if one considers the support of $\phi^*$ and $\hat\phi_n$. Thus, we continue with sign consistency.

We shall need a few more definitions. Let $s=|S|$ denote the cardinality of the set of true nonzero coefficients, and let $\nu=p-s=|S^c|$, with $S^c = \{1,\ldots,p\} \setminus S$. For a set of indexes $I$, we will write $x_I = \{ x_i, i \in I \}$ to denote the subvector of $x$ whose elements are indexed by the coordinates in $I$. Similarly, $x_Iy_I$ is a vector with elements $x_iy_i$. %and $\|x_I\|$ is (some) norm of $x$ restricted to $I$, i.e., $\|x_I\|_\infty$ is $\max_{i\in I} |x_i|$. 
For a $n\times p$ design matrix $X$, we let $X_I$, for any subset $I$ of $\{1,2,\ldots, p\}$, denote the sub-matrix of $X$ with columns as indicated by $I$. Sub-matrices of the autocovariance matrix $\Gamma_p$ (and of any other matrix), are denoted similarly. For example, $\Gamma_{II^c}$ is $(\gamma(i-j))_{i\in I, j\notin I}$.  Finally, let $\alpha_n=\min_{j\in S}|\phi_j^*|$ denote the magnitude of the smallest nonzero coefficient.
Finally,  although virtually all quantities related to (\ref{eq:lasso_ar}) depend on $n$, we do not always make this dependence explicit in our notation.

%Some quantities, say $x\equiv x_n$ may depend on $n$, the sample size. In this case we write $x_{n,I}$. 

We are now ready to present our first result:

\begin{thm}\label{thm:sign}
Consider the settings of the AR($p$) process describe above. Assume that
\begin{description}
\item[(i)] there exists a finite, positive constant $C_{\max}$ such that $\|\Gamma^{-1}_{SS}\|\leq C_{\text{max}}$;
\item[(ii)] there exists an $\epsilon\in(0,1]$ such that $\|\Gamma_{S^cS}\Gamma^{-1}_{SS}\|_\infty\leq 1-\epsilon$.
\end{description}
Further, assume that the following conditions hold:
\begin{equation}\label{eq:sign_cond1}
\limsup_{n\rightarrow \infty} \frac{\max_{i\in S}\lambda_{n,i}}{\min_{j\in S^c}\lambda_{n,j}} \leq 1 \; ,
\end{equation}
\begin{equation}\label{eq:sign_cond2}
\frac{1}{\alpha_n} \Big[\sqrt{s/n} +
 \lambda_n\|\lambda_{n,S}\|_\infty\Big] \longrightarrow 0 \quad ,\quad \text{as} \quad n\rightarrow\infty \; ,
\end{equation}
\begin{equation}\label{eq:sign_cond3}
\frac{n\lambda_n ^2 (\min_{i\in S^c} \lambda_{n,i})^2}{\max\{s,\nu\}}\longrightarrow \infty \qquad , \qquad \text{as}\quad  n\rightarrow \infty \; .
\end{equation}
Let $p=\OO{\log n}$.
Then, the Lasso estimator $\hat\phi_n$ is sign consistent (cf. (\ref{eq:sign_consistency})). 
\end{thm}

Condition {\bf (ii)} in Theorem \ref{thm:sign} is assumed in various guises elsewhere in the Lasso literature (see, e.g., \cite{wainwright}, \cite{zhaobin} and \cite{adaptive}). It is an incoherence condition, which controls the amount of correlation between relevant variables and irrelevant variables.  %\citet{zhaobin} use a similar condition which they call the irrepresentable condition. 
Condition (\ref{eq:sign_cond1}) is intuitively clear and it appears under similar form in \cite{NardiRinaldo}. It captures the rationale , recalling that one may have $\lambda_j<\lambda_k $ for $j<k$, that (even) the largest penalty coefficient of the relevant lags should be kept asymptotically smaller than the smallest penalty coefficient of the irrelevant lags. Conditions (\ref{eq:sign_cond2}) and (\ref{eq:sign_cond3}) are similar to conditions appearing in \cite{wainwright}, \cite{NardiRinaldo}, and \cite{LaffertySpam}, to name but a few. The fraction $\sqrt{s/ n}$ in (\ref{eq:sign_cond2}) is in line with similar works, mentioned above. For example, under the \emph{linear sparsity scheme}, i.e., $s=\alpha  p$, with $\alpha\in (0,1)$ (see \cite{wainwright}), and with $p$ comparable to $n$, the Gaussian ensemble leads to a fraction of order $\OO{\log n/n}$, which is similar to the fraction under the current setting, for which we have $p=\OO{\log n}$.

\subsection{Estimation and Prediction Consistency}\label{sec:estimation_prediction}

Our next result is about \emph{estimation consistency}. An estimator $\hat \phi_n$ is said to be estimation consistent, or $l_2$-consistent if $\|\hat\phi_n -\phi^*\|$ converges to zero, as $n$ tends to infinity. We have the following:
\begin{thm}\label{thm:estimation}
Recall the settings of the AR($p$) process set forth below (\ref{eq:ar_model}). Let $p=\OO{\log n}$, and $\alpha_n=p^{1/2}(n^{-1/2}+\lambda_n\|\lambda_{n,S}\|)$. Assume that $\lambda_n\|\lambda_{n,S}\|=\OO{n^{-1/2}}$. Then, the Lasso estimator $\hat\phi_n$ is estimation consistent with a rate of order $\OO{\alpha_n}$.
\end{thm}

\emph{Prediction consistency} is about a similar convergence statement, but for the prediction of future values using the fitted model. Formally, prediction consistency holds if $\|X\hat\phi_n-X\phi^*\|$ converges to zero, as $n$ tends to infinity. We show below a similar result when the sample autocovariance matrix $X'X$ is replaced by the (theoretical) autocovariance matrix $\Gamma_p$. The autoregressive settings assumed here are, in some sense, much more challenging than in linear (parametric or non-parametric) regression models, for two reasons. Firstly, the design matrix is not fixed as is usually assumed, and  secondly, the entries of the $X$  are not independent across rows, as is usually assumed for random designs.

The family of AR processes considered here are, in fact, a subset of a larger family of time series. In order to establish the prediction consistency result, we make an explicit use of the structure  of this larger family, to which we now describe. 

Following \cite{zeevi}, we denote by $\mathcal{H}_\rho (l,L)$, for some $\rho>1$, $0<l<1$, and $L>1$, a family consisting of all stationary Gaussian time series with $\mathbb{E}X_t=0$, $\mathbb{E}|X_t|^2=1$, and with
\[
0<l\leq |\psi(z)|\leq L \; ,
\]
for every complex $z$ with $|z|\leq \rho$, where $\psi(z)$ is the MA($\infty$) transfer function related to the AR polynomial by $\psi(z)=1/\phi(z)$.

We shall need the notion of a strong mixing (or $\alpha$-mixing) condition. Let $\{X_t\}$ be a time series defined on a probability space $(\Omega, \mathcal{F}, \mathbb{P})$. For any two (sub) $\sigma$-fields $\mathcal{A}$ and $\mathcal{B}$, define 
\[
\alpha(\mathcal{A},\mathcal{B})=\sup_{A\in\mathcal{A}, B\in\mathcal{B}} |\mathbb{P}(A\cap B)-\mathbb{P}(A)\mathbb{P}(B)| \; .
\]
Denote by $\mathcal{F}_s^t$, the $\sigma$-field generated by $(X_s, \ldots, X_t)$, for $-\infty\leq s\leq t\leq\infty$. Then, $\{X_t\}$ is said to be strongly mixing if $\alpha_X(m)\rightarrow 0$, as $m\rightarrow\infty$, where
\[
\alpha_X(m) = \sup_{j\in\{0,\pm 1,\pm 2, \ldots\}} \alpha(\mathcal{F}_{-\infty}^j,\mathcal{F}_{j+m}^\infty) \; .
\]

Attractiveness of $\mathcal{H}_\rho (l,L)$ comes from the fact that processes in $\mathcal{H}_\rho(l,L)$ are strong mixing with an exponential decay, i.e.
\begin{equation}\label{eq:decay_alpha}
\alpha_X(m) \leq 2\left(\frac{L\rho}{l(\rho-1)}\right)^2 \rho^{-m} \; .
\end{equation}
This follows since processes in $\mathcal{H}_\rho(l,L)$ have exponentially decaying AR coefficients as well as exponentially decaying autocovariances (see \cite[Lemma 1, and in particular, expression (39)]{zeevi}).

For every $p$-dimensional vector $a$ and $p\times p$ symmetric matrix $A$, we denote with  $\|a\|_A^2=a'Aa$, the (squared) $l_2$-norm associated with $A$. Let $C_1, C_2$ be two universal constants (their explicit values are given within the proof of the following theorem). Define
\begin{equation}\label{eq:constants}
\beta_1 = 1+\frac{1}{\log\rho} \quad , \quad \beta_2 = 1+\frac{L\rho}{l(\rho-1)} \; , \qquad \text{and} \quad D=(C_1^3C_2\beta_1^2\beta_2^3)^{1/5} \; .
\end{equation}
Let $\lambda_{\min}=\min_{j=1,\ldots, p} \lambda_{n,j}$, and $\lambda_{\max}=\max_{j=1,\ldots, p} \lambda_{n,j}$. We have:
\begin{thm}\label{thm:prediction}
Recall the settings of the AR($p$) process set forth below (\ref{eq:ar_model}). Let $p=\OO{\log n}$. Assume:
\begin{description}
\item[(i)] There exists a finite, positive constant $M$ such that $\lambda_{n,j}\leq M$, for every $j=1,\ldots, p$.
\item[(ii)] For every $p\geq 2$, there exists a positive constant $\kappa_p$, such that
\[
\Gamma_p-\kappa_p \text{diag}(\Gamma_p)
\]
is a positive semi-definite matrix.
\end{description}
If $\lambda_n (s/p)^{1/2}\leq D n^{-2/5}$, then there exist a constant $C$ (depending only on $M$), and constants $F_1$ and $F_2$ (depending only on $C_1,C_2,\beta_1,\beta_2$), such that for all $0<c<\infty$, and all $y>\sigma^2(n+Dn^{3/5})$,
\[
\|\hat\phi_n-\phi^*\|_{\Gamma_p}^2\leq C\lambda_n^2\frac{s}{\kappa_p}
\] 
holds true with probability at least $1-\pi_n$, where
\begin{equation}\label{eq:pi}
\pi_n \leq 6p\exp\left\{ -F_1\min \left\{(\sigma^{-2}y-n)^{1/3}, c^2\sigma^{-2}, \frac{n^2\lambda_n^2\lambda^2_{\min}}{y+cn\lambda_n\lambda_{\max}/2}\right\}  \right\} + p^2\exp\left\{ -F_2 n\lambda_n^2 (s/p^2)\right\} \; .
\end{equation}
\end{thm}

Condition {\bf (ii)} has been used in the context of aggregation procedures for nonparametric regression with fixed design (\cite{buneaFixed}), and also for nonparametric regression with random design (\cite{buneaRandom}).

Theorem \ref{thm:prediction} may be utilized to show that the Lasso estimator $\hat\phi_n$ is prediction consistent. One only needs to make sure that the decay of the bound (\ref{eq:pi}) on $\pi_n$. The theorem actually gives a whole range of possible rates of decay, by picking $c$ and $y$. One \emph{possible} choice is given below.
\begin{cor}\label{cor:prediction}
Let $\lambda_n = n^{-\alpha}$, with $\alpha\in(2/5,1/2)$. Let $c=D_1y/(n\lambda_n\lambda_{\max})$, and $y=D_2n$, for positive constants $D_1, D_2$. If $(s/p)^{1/2}\leq D n^{\alpha-2/5}$, then there exists an appropriate constant $F$, such that the bound (\ref{eq:pi}) on $\pi_n$ is smaller than
\[
p^2\exp \Big\{ -F\min\Big\{ n^{1/3},n^{2\alpha}/\lambda^2_{\max}, n^{1-2\alpha}\lambda_{\min}^2, n^{1-2\alpha}s/p^2 \Big\}   \Big\} \; ,
\]
which tends to zero as $n$ goes to infinity.
\end{cor}

\section{Illustrative Simulations}\label{sec:simul}

We consider a sparse autoregressive time series of length 1000 obeying the model
\begin{equation}\label{eq:simu}
X_t = 0.2 X_{t-1} + 0.1 X_{t-3} + 0.2 X_{t-5} + 0.3 X_{t-10} +  0.1 X_{t-15} + Z_t,
\end{equation}
with nonzero coefficients at lags 1, 3, 5, 10 and 15, where the innovations $Z_t$  are i.i.d. Gaussians with mean zero and standard deviation $0.1$. The coefficients were chosen to satisfy the characteristic equation for a stationary AR process.

\begin{figure}[!ht]
	\centering
	\includegraphics[width=4in,height=7in]{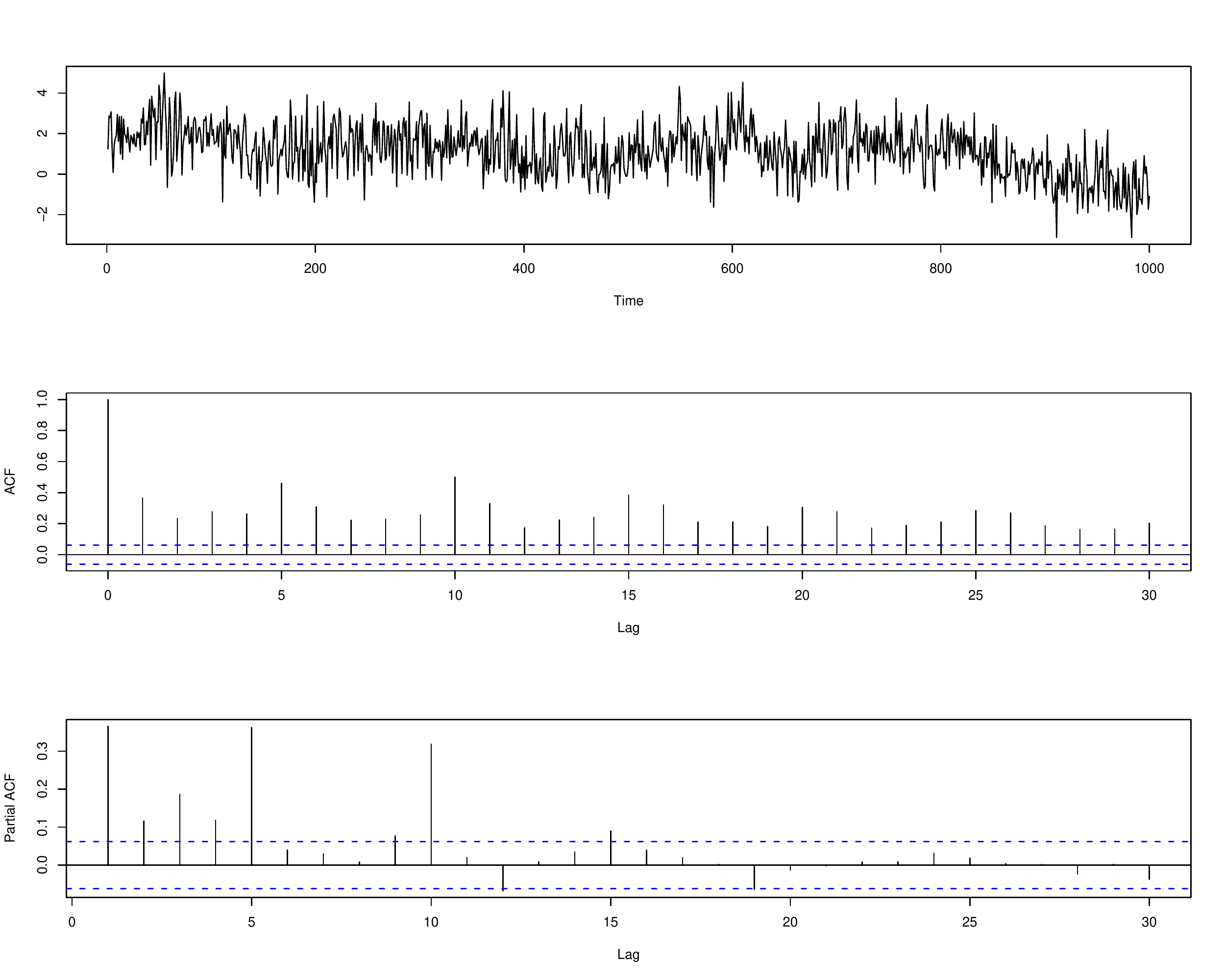}
	\caption{A time series simulated from the sparse autoregressive model (\ref{eq:simu}) along with its autocorrelation and partial autocorrelation coefficients.}
	\label{fig:series}
\end{figure}

\begin{figure}[t]
	\centering
	\includegraphics[width=4in,height=2in]{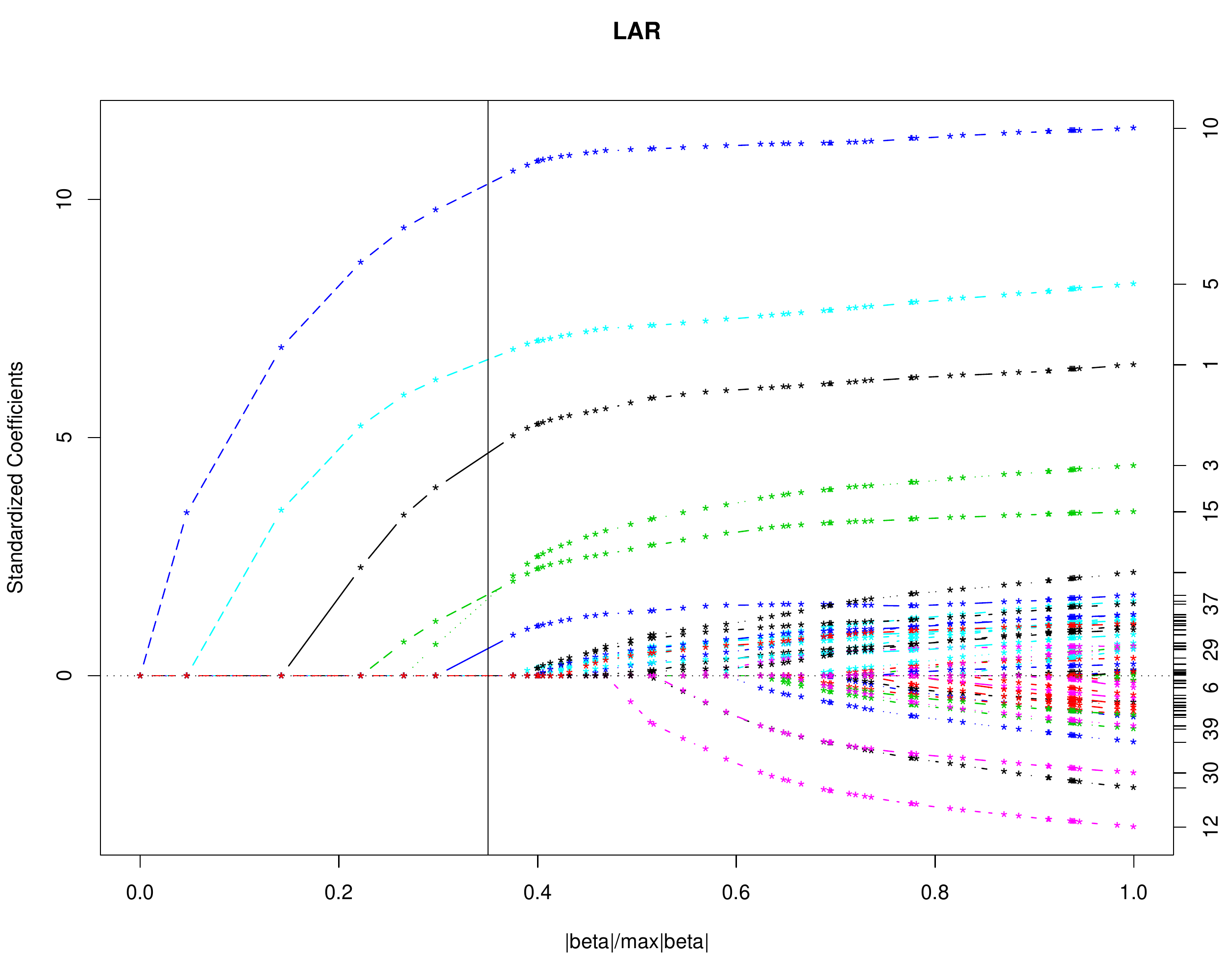}
	\caption{Solution paths of the {\tt lars} algorithm when applied to the time series displayed in Figure \ref{fig:series}. The vertical bar represents the optimal $\ell_1$ penalty for this time series selected using cross validation.}
	\label{fig:lars}
\end{figure}

Figure \ref{fig:series} shows one time series simulated according to the model (\ref{eq:simu}), along with its autocorrelation and partial autocorrelation plots.
For this time series, Figure \ref{fig:lars} shows the solution paths computed using the {\tt R} algorithm {\tt lars} and for a value of $p=50$.  Notice that we only use one penalty parameter, i.e. we penalize equally all the autoregressive coefficients. The vertical line marks the optimal $\ell_1$ threshold found by cross validation. In our simulations, we declared significant the variables whose coefficients have nonzero solution paths meeting the vertical line corresponding to the cross validation value. 

\begin{figure}[!ht]
	\centering
	\includegraphics[width=4in,height=2in]{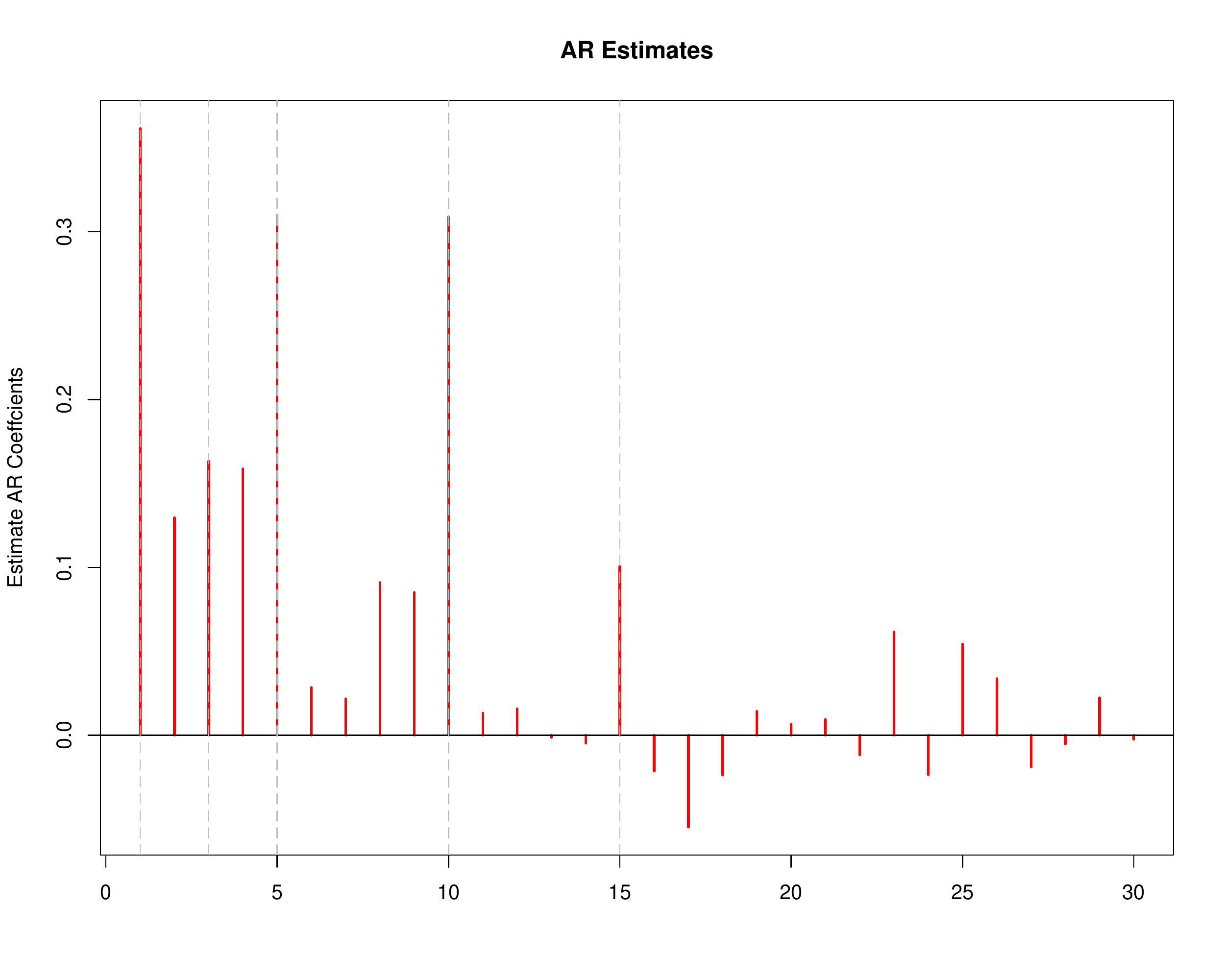}
	\caption{Autoregressive coefficients for the time series of Figure \ref{fig:series} obtained using the routine {\tt ar}. The dashed vertical line mark the lags for true nonzero coefficients.}
	\label{fig:ar}
\end{figure}

Notice that, in the exemplary instance displayed in Figure \ref{fig:lars}, all the nonzero autoregressive coefficients are correctly included in the model. Furthermore, a more careful inspection of the solution paths reveals that the order at which the significant variables enter the set of active solutions match very closely the magnitude of the coefficients used in our model, with $\phi_{10}$ and $\phi_5$, the more significant coefficients, entering almost immediately, and $\phi_3$ and $\phi_{15}$ entering last.
In contrast, Figure \ref{fig:ar} displays the fitted values for the first 30 autoregressive coefficients computed  using the Yule-Walker method implemented using {\tt R} by the routine {\tt ar} (note that the Yule-Walker estimator has the same asymptotic distribution as the MLE's). Notice that the solution is non-sparse. The dashed vertical gray lies indicate the true nonzero coefficients. The autoregressive order of the model was correctly estimated to be 15 using the AIC criterion.

\begin{figure}[!ht]
	\centering
	\begin{tabular}{cc}
	\includegraphics[width=3in,height=2in]{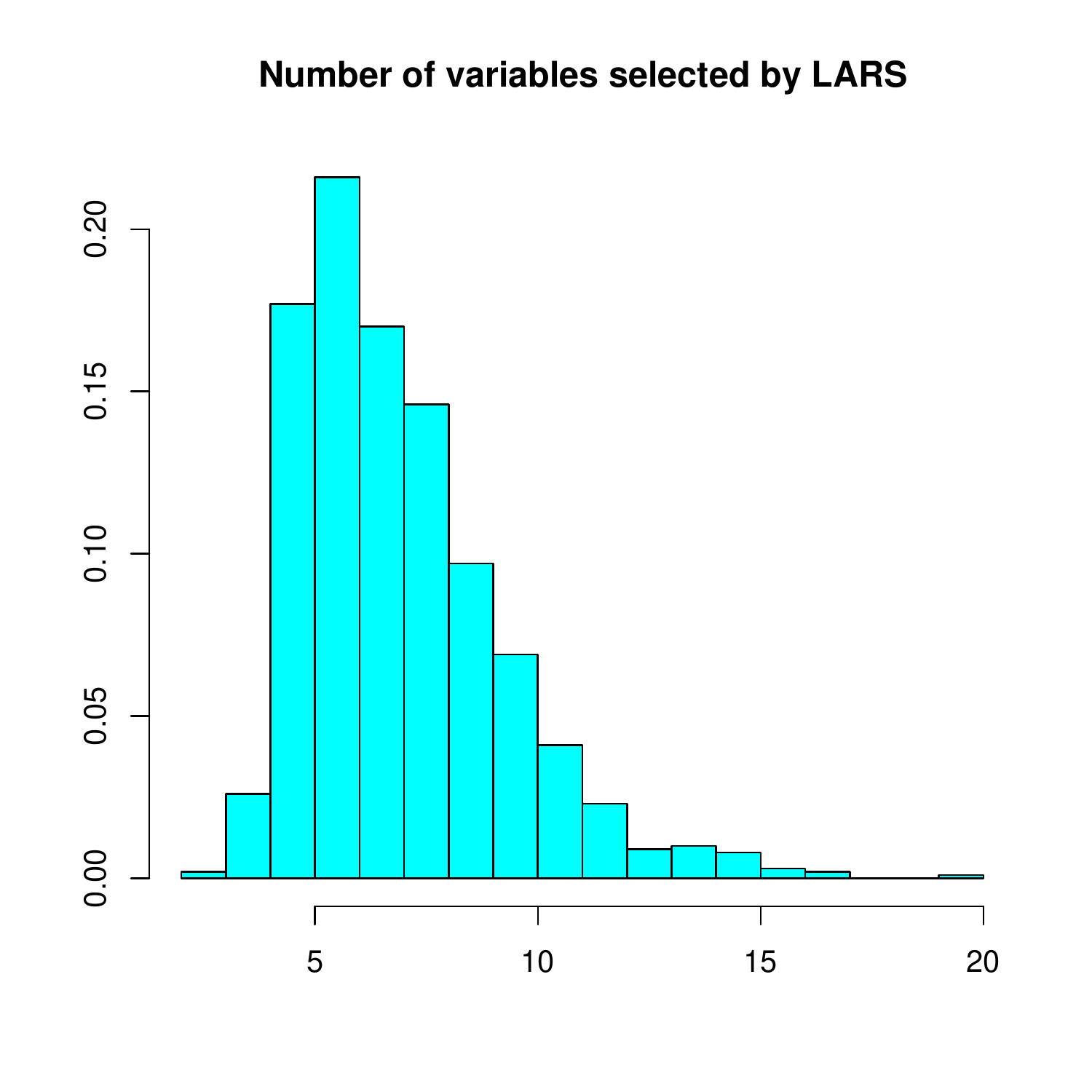} & \includegraphics[width=3in,height=2in]{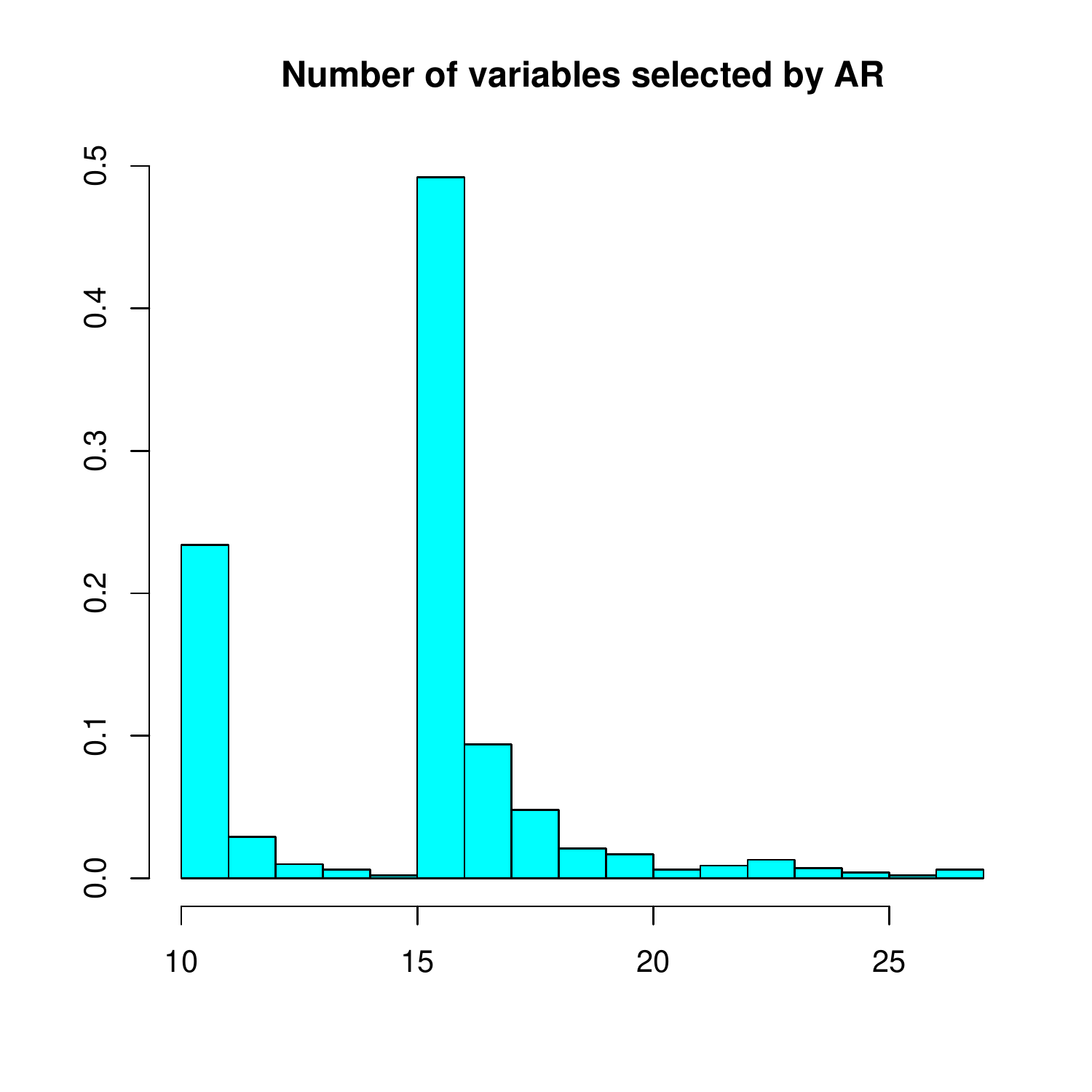}\\
	{\bf a)} & {\bf b)}\\
	\end{tabular}
	\caption{Distribution for the number of variables selected by {\bf a)} the {\tt lars} algorithm using cross validation and {\bf b)} the {\tt ar} algorithm using AIC over 1000 simulations of the time series described in (\ref{eq:simu}).}
	\label{fig:num}
\end{figure}

We simulated 1000 time series from the model (\ref{eq:simu}) and we selected the significant variables according to the cross-validation rule described above. Figure \ref{fig:num} {\bf a)} displays the histogram of the number of selected variables. The mean and standard deviations of these numbers are 6.42 and 2.44, respectively, while the minimum, median and maximum numbers are 3, 6 and 22, respectively. In comparison, Figure \ref{fig:num} {\bf b)} shows the histogram of the autoregressive orders determined by AIC in {\tt ar}.
Table \ref{tab:simu} displays some summary statistics of our simulations. In particular, the second row shows the number of times, out 1000 simulated time series, that each of the nonzero autoregressive coefficients was correctly selected. The second row indicates the number of times the variable corresponding to each nonzero  coefficient in (\ref{eq:simu}) was among the first five selected variables. Notice that  $\phi_{10}$ and $\phi_5$ are always included among the selected variables, while $\phi_3$ and $\phi_{15}$ have a significantly smaller, but nonetheless quite high, chance of being selected.

\begin{table}
\centering
\caption{Number of times the nonzero autoregressive coefficients are correctly identified and number of times they are correctly selected among the first $5$ variables entering the solution paths.}
\label{tab:simu}
\begin{tabular}{cccccc}
\hline
$\phi$ & $\phi_1$ & $\phi_3$ & $\phi_5$ & $\phi_{10}$ & $\phi_{15}$\\
\hline
\multicolumn{1}{c}{Value} & 0.2 & 0.1 & 0.2 & 0.3 & 0.1\\ 
\multicolumn{1}{c}{Number of times correctly selected} &  992 & 754 & 1000 & 1000 & 913\\
\multicolumn{1}{c}{Number times selected among first 5} &  992 & 602 & 1000 & 1000 & 895\\
\hline
\end{tabular}
%\end{center}
\end{table}

%\begin{table}\label{tab:simu}
%\begin{center}
%\begin{tabular}{c|ccccc|}
%%& \multicolumn{5}{|c|}{Coefficients} \\
%& $\phi_1$ & $\phi_3$ & $\phi_5$ & $\phi_{10}$ & $\phi_{15}$\\
%\hline
%\multicolumn{1}{|c|}{Value} & 0.2 & 0.1 & 0.2 & 0.3 & 0.1\\ 
%\hline
%\multicolumn{1}{|c|}{Number of times correctly selected} &  992 & 754 & 1000 & 1000 & 913\\
%\hline
%\multicolumn{1}{|c|}{Number times selected among first 5} &  992 & 602 & 1000 & 1000 & 895\\
%\hline
%\end{tabular}
%\end{center}
%\caption{Number of times the nonzero autoregressive coefficients are correctly identified and number of times they are correctly selected among the first $5$ variables entering the solution paths.}
%\end{table}

\begin{figure}[!h]
\centering
\begin{tabular}{cc}
\includegraphics[width=2in]{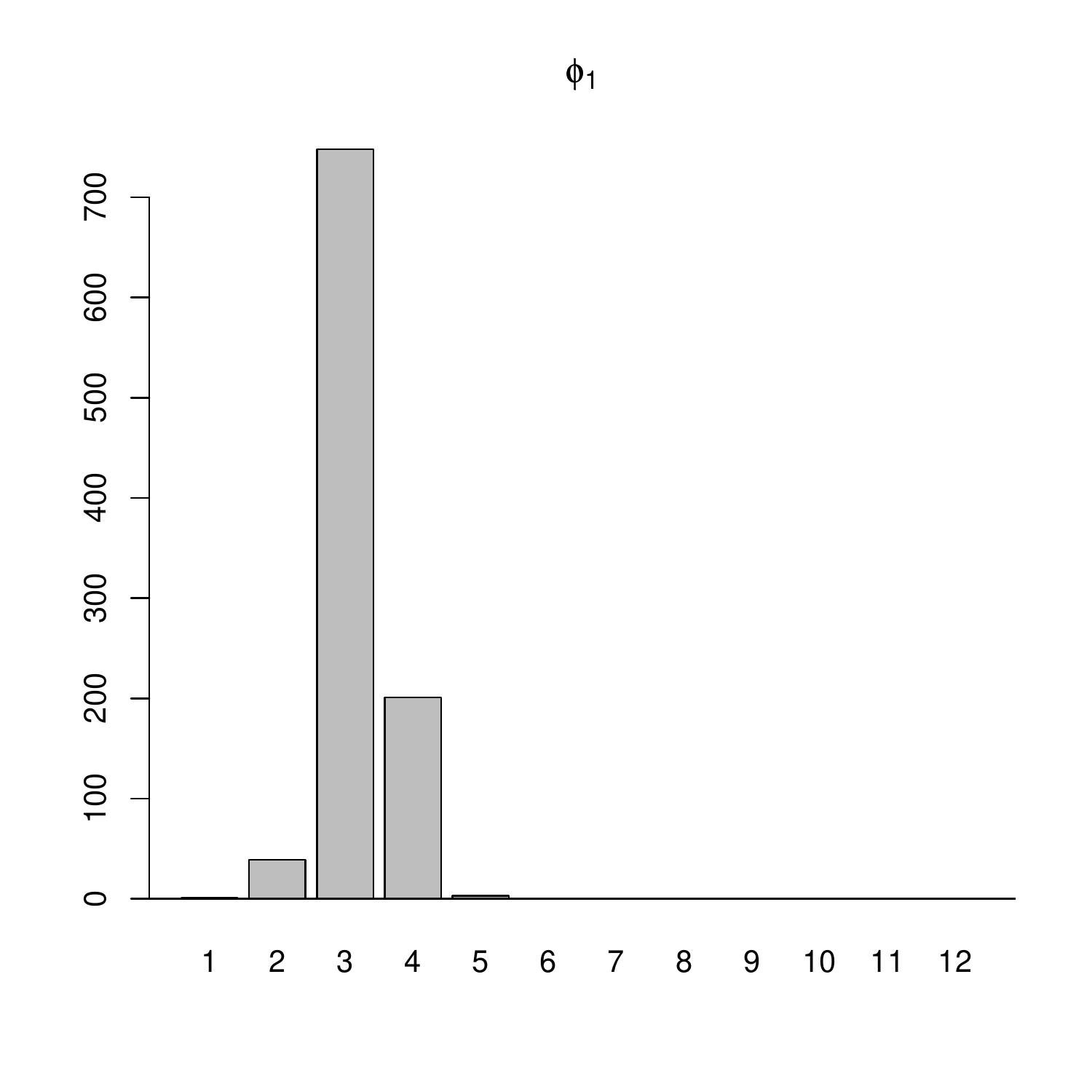} &  \includegraphics[width=2in]{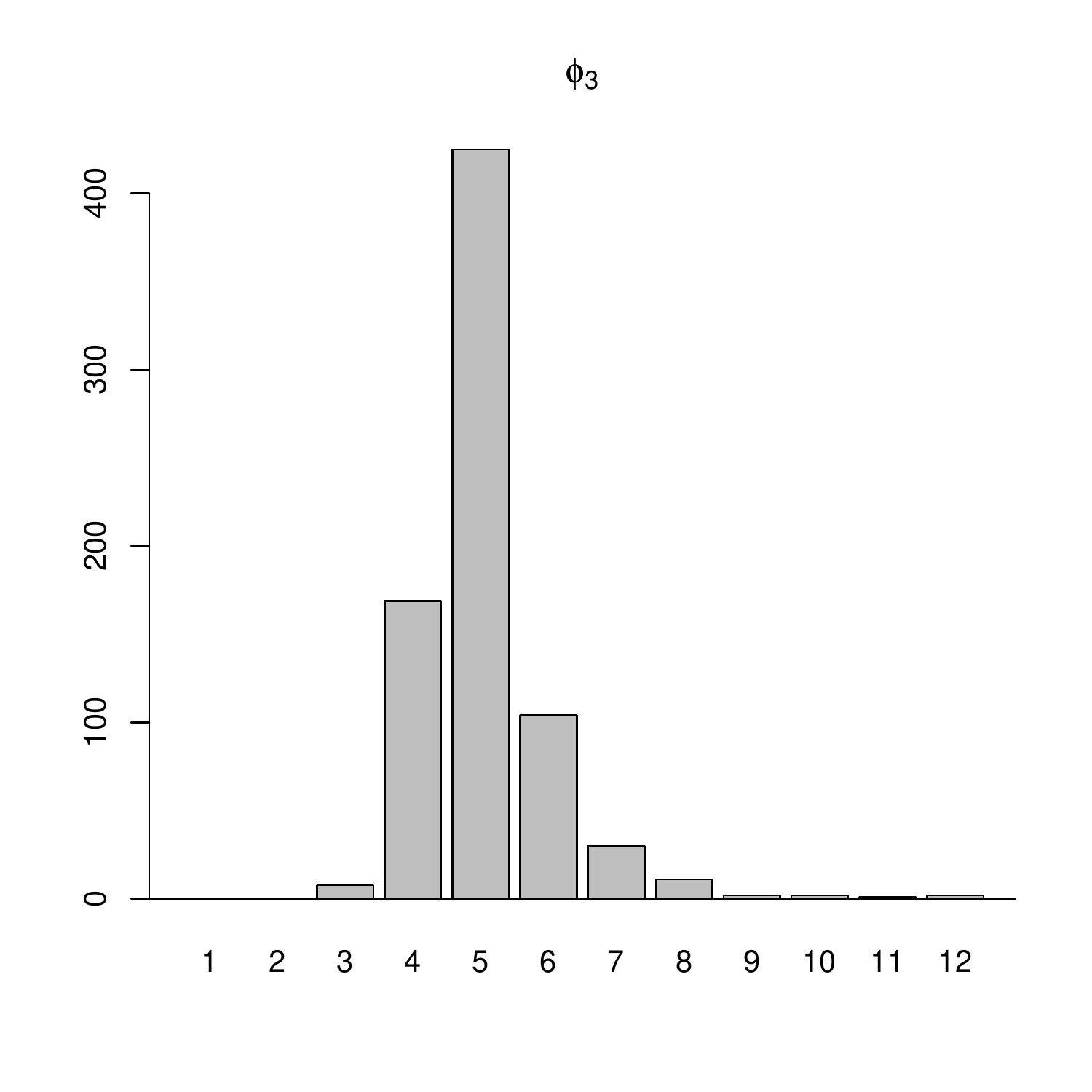} \\
\includegraphics[width=2in]{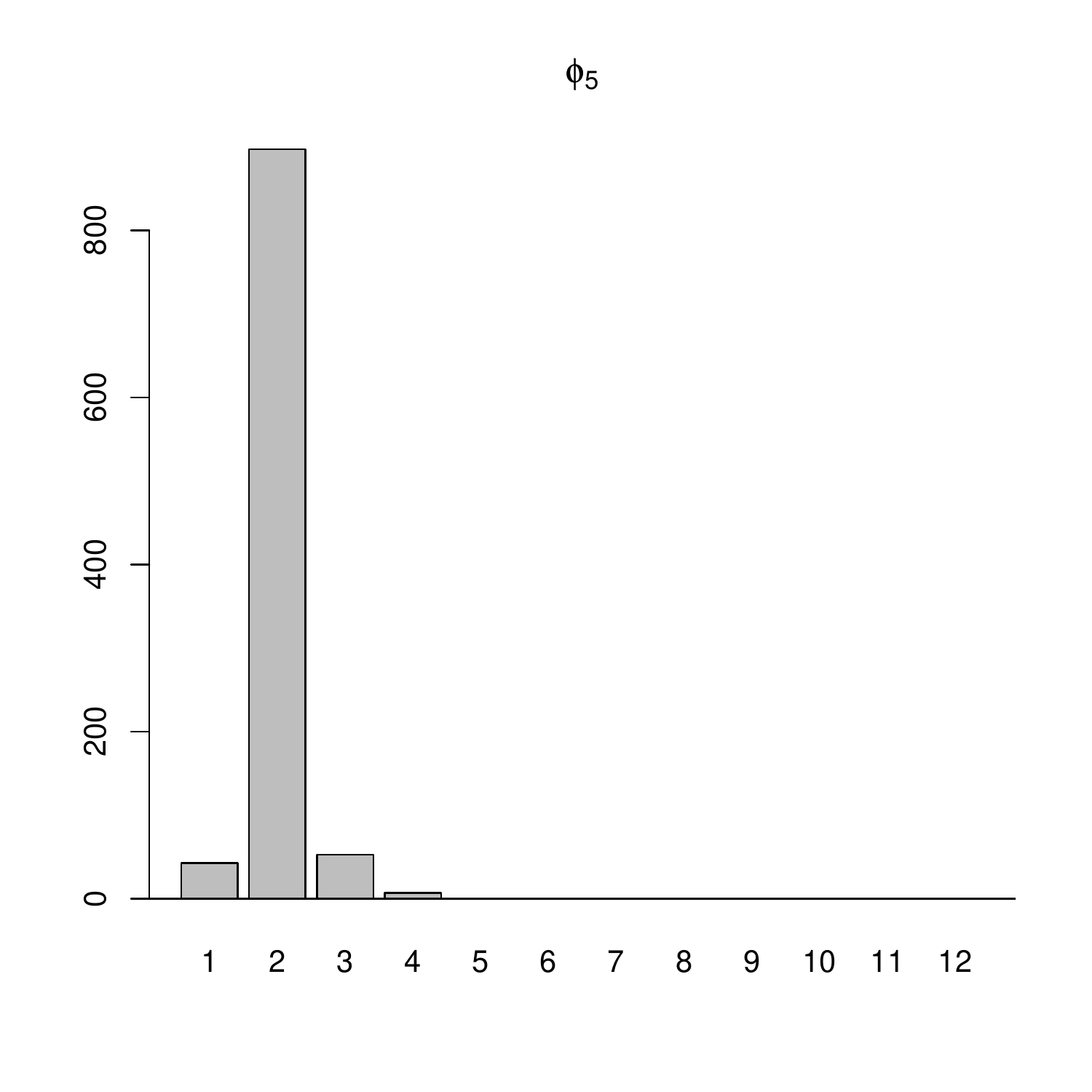} &  \includegraphics[width=2in]{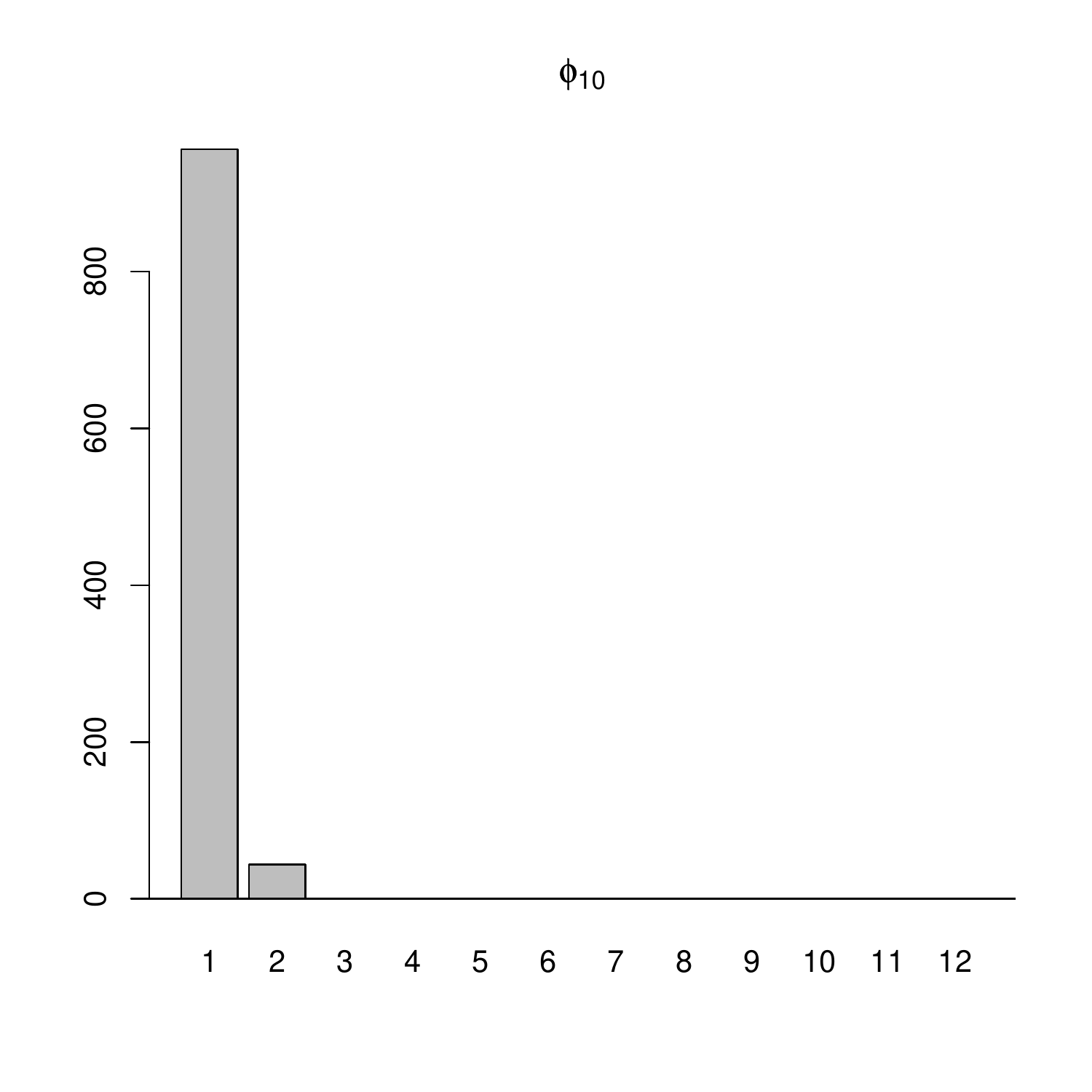}  \\
\includegraphics[width=2in]{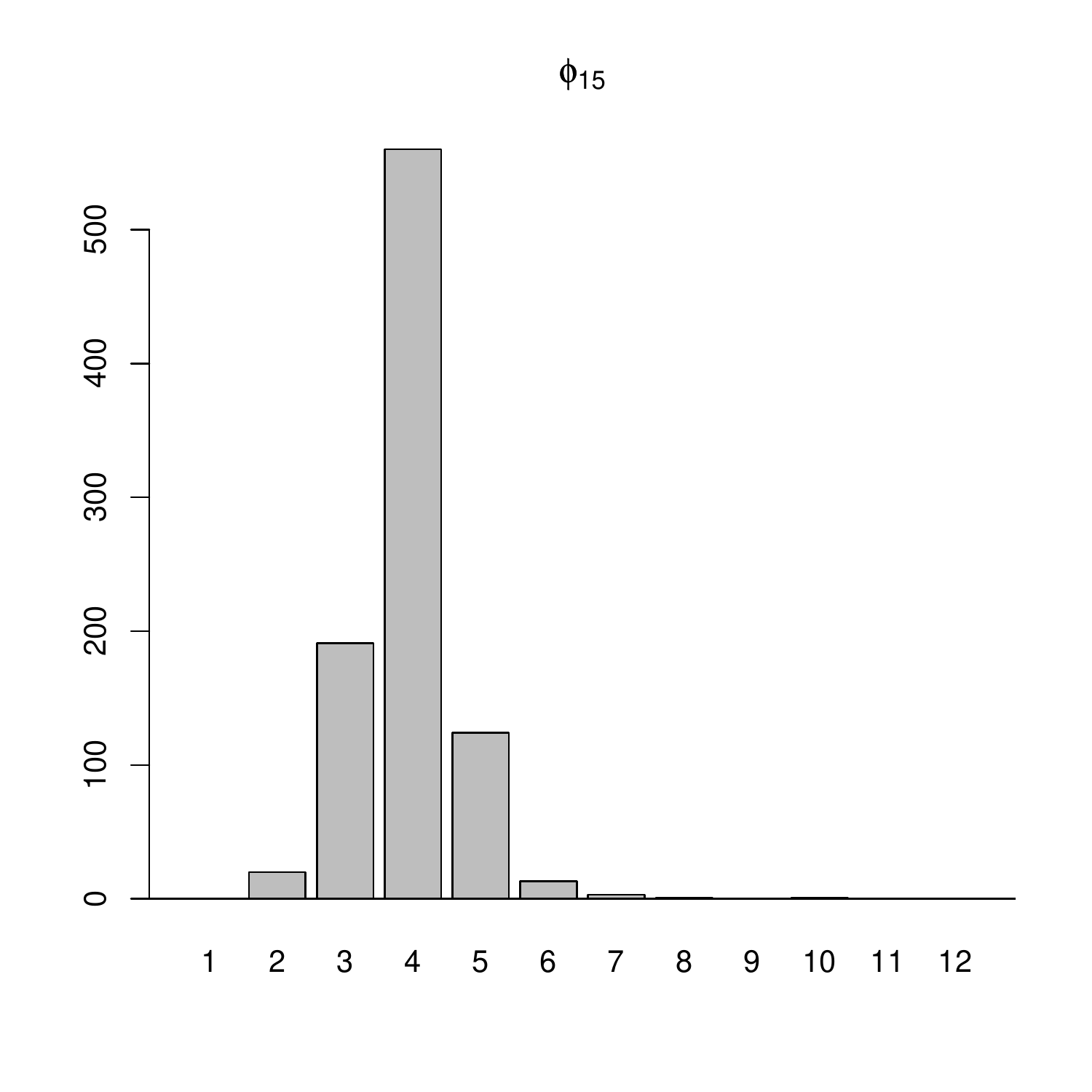} &   \\
\end{tabular}
\caption{Frequencies of the order at which the 5 autoregressive coefficients entered the solutions paths for the {\tt lars} algorithm over 1000 simulations of the time series described in (\ref{eq:simu}).}
\label{fig:phi}
\end{figure}

We also investigated the order at which the autoregressive coefficients entered the solution paths, the rationale being that more significant nonzero variables enter sooner, in accordance with the way the {\tt lars} algorithm works (see \cite{lars}). Figure \ref{fig:phi} summarizes our findings. In each of the barplots, the $x$-axis indexes the steps at which the variable corresponding to the autoregressive coefficient enters the solution path, while the $y$-axis displays the frequency. Interestingly enough, in most cases, $\phi_{10}$ and $\phi_{5}$ are selected as the first and second nonzero variables, while $\phi_{15}$ and, in particular, $\phi_3$ enter the set of active variables later and are not even among the first five variables selected in  1.9\% the and 20.2\%  of cases, respectively.

%\multicolumn{3}{|c||}{} 

%\begin{figure}[!ht]
%	\centering
%	\includegraphics[width=6in,height=3in]{lasso2.pdf}
%	\caption{blablabla.}
%	\label{fig:lasso2}
%\end{figure}

\section{Discussion}\label{sec:discussion}

We defined the Lasso procedure for fitting an autoregressive model, where the maximal lag may increase with the sample size. Under this double asymptotic framework, the Lasso estimator was shown to possess several consistency properties. In particular, when $p=\OO{\log n}$, the Lasso estimator is model selection consistent, estimation consistent, and prediction consistency. The advantage of using the Lasso procedure in conjunction with a growing $p$ is that the fitted model will be chosen among all possible AR models whose maximal lag is between $1$ and $\OO{\log n}$. Letting $n$ go to infinity, we may virtually obtain a good approximation for a general linear time series.

As mentioned in Section \ref{sec:ar}, the assumption about Gaussianity of the underlying noise $\{Z_t\}$ is not necessary. The proof of the model selection consistency result (Theorem \ref{thm:sign}) avoids making use of Gaussianity by using Burkholder's inequality in conjunction with a maximal moment inequality. The proof of the estimation consistency result (Theorem \ref{thm:estimation}) requires Lemma \ref{lem:consis_X}, which does make use of the assumed Gaussianity. However, this is not crucial. In fact, we can relax the Gaussianity assumption and require only the $Z_t$ are $ IID(0,\sigma^2)$, with bounded fourth moment (see~\cite[p. 226-227]{MR1093459}). In this case, instead of using Wick's formula we may apply the moving average representation $X_t=\sum_{j=0}^\infty \psi_j Z_{t-j}$, along with the absolute summability of the $\psi_j$'s. Finally, the prediction consistency result (Theorem \ref{thm:prediction} and Corollary \ref{cor:prediction}) may also be obtained by relaxing the Gaussianity assumption. One only needs to impose appropriate moment conditions of the driving noise.

The autoregressive modeling via the Lasso procedure stimulates other interesting future directions. In many cases, non-linearity is evident from the data. In order to capture deviation from linearity, one may try to fit  a non-linear (autoregressive) time series model to the data in the form
\[
X_t = \phi_1X_{t-1}+\cdots +\phi_pX_{t-p}+\sum_{\nu=2}^p\{\phi^{i_1,\ldots , i_\nu}\prod_{j=1}^\nu X_{t-i_j}\}+Z_t \; ,
\]
where we used the Einstein notation for the term in the curly brackets, to indicate summation over all $i_1<i_2<\ldots <i_\nu$. Notice that for even mild values of $p$, the number of possible interaction terms may be very large. This is a very challenging problem as one needs to obtain a solid understanding of the properties of the non-linear autoregressive process before applying the Lasso (or any other) procedure.

\section{Proofs}  \label{sec:proofs}

Here we prove Theorems \ref{thm:sign}, Theorem \ref{thm:estimation}, and Theorem \ref{thm:prediction}. Recall scheme (\ref{eq:lasso_ar}). This is a convex minimization problem. Denote by $M_{\Lambda_n}(\cdot)$, for $\Lambda_n=\{\lambda_n, \{\lambda_{n,j}\, ,\, j=1,\ldots p\}\}$, the objective function, i.e.,
\begin{equation}\label{eq:objective_function}
M_{\Lambda_n}(\phi) = \frac{1}{2n}\|y-X\phi\|^2
+\lambda_n \sum_{j=1}^p \lambda_{n,j}|\phi_j| \; .
\end{equation}
The Lasso estimator is an optimal solution to the problem  $\min\{M_{\Lambda_n}(\phi) \, , \, \phi\in\mathbb{R}^p\}$. Gradient and Hessian of the least-squares part in $M_{\Lambda_n}(\cdot)$ are given, respectively, by $n^{-1}\mathfrak{X}\phi-n^{-1}\sum_{t=1}^nX_t{\bf X_t}$, and $n^{-1}\mathfrak{X}$, where $\mathfrak{X}$ (the gram matrix associated with the design matrix $X$), and ${\bf X_t}$ is a notation that we use throughout this section:
\[
\mathfrak{X} = X'X \quad , \quad {\bf X_t} = (X_{t-1}, \ldots, X_{t-p})' \; .
\]

\subsection{Model Selection Consistency}

\begin{proof}[Proof of Theorem \ref{thm:sign}]
We adapt a Gaussian ensemble argument, given in \cite{wainwright}, to the present setting. Standard optimality conditions for convex optimization problems imply that $\hat\phi_n\in\mathbb{R}^p$ is an optimal solution to the problem $\min\{M_{\Lambda_n}(\phi) \, , \, \phi\in\mathbb{R}^p\}$, if, and only if,
\begin{equation}\label{eq:optim_condin}
\frac{1}{n}\mathfrak{X}\hat\phi_n -\frac{1}{n}\sum_{t=1}^n {X_t\bf X_t }+\lambda_n\hat\xi_n
 = 0 \; ,
\end{equation}
where $\hat\xi_n\in\mathbb{R}^p$ is a sub-gradient vector with elements $\hat\xi_{n,j}=
\text{sgn}(\hat\phi_{n,j})\lambda_{n,j} $ if $\hat\phi_{n,j}\ne0$, and $|\hat\xi_{n,j} |\leq \lambda_{n,j}$ otherwise. Plugging the model
structure, $y=X\phi^*+Z$, into~(\ref{eq:optim_condin}), one can see that the optimality
conditions become
\begin{equation}\label{eq:optim_condin2}
\frac{1}{n}\mathfrak{X}(\hat\phi_n -\phi^*)-\frac{1}{n}\sum_{t=1}^n {Z_t\bf X_t }+\lambda_n\hat\xi_n
 = 0 \; .
\end{equation}

Recall the sparsity set, $S=\{j  \, : \, \phi_j^*\ne 0\}=\text{supp}(\phi^*)$, the sparsity cardinality $s=|S|$, and $\nu=p-s=|S^c|$. Decomposing the design matrix $X$ to relevant and non-relevant variables, $X=(X_S , X_{S^c})$, we may write the gram matrix $\mathfrak{X}$ as a block matrix of the form
\[
\mathfrak{X} = \left( \begin{array}{cc}
\mathfrak{X}_{SS} & \mathfrak{X}_{SS^c} \\ \mathfrak{X}_{S^cS} & \mathfrak{X}_{S^cS^c} \\ \end{array}\right)  =  \left( \begin{array}{cc}
X_S'X_S & X_S'X_{S^c} \\ X_{S^c}'X_S & X_{S^c}'X_{S^c} \\ \end{array}\right) \; .
\]
Notice, for example,  that  $\mathfrak{X}_{SS}=\left(\sum_{t=1}^n X_{t-i}X_{t-j}\right)_{i,j\in S }$. Incorporating this into the optimality conditions (\ref{eq:optim_condin2}) we obtain the following two relations,
\begin{eqnarray*}
\frac{1}{n}\mathfrak{X}_{SS}[\hat\phi_{n,S}-\phi^*_S]-\frac{1}{n}\sum_{t=1}^nZ_t{\bf X_t^S} &=& -\lambda_n\lambda_{n,S}\,\text{sgn}(\phi_S^*)  \; ,\\
\frac{1}{n}\mathfrak{X}_{S^cS}[\hat\phi_{n,S}-\phi_S^*]-\frac{1}{n}\sum_{t=1}^nZ_t{\bf X_t^{S^c}} &=& -\lambda_n\hat\xi_{n,S^c} \; ,
\end{eqnarray*}
where ${\bf X_t^{S}}$, and ${\bf X_t^{S^c}}$ are vectors with elements \{$X_{t-i}\, , \, i\in S\}$, and \{$X_{t-i}\, , \, i\in S^c\}$, respectively. 
If $n-s\geq s$ then $\mathfrak{X}_{SS}$ is non-singular with probability one, and we can solve for $\hat\phi_{n,S}$ and $\hat\xi_{n,S}$,
\begin{eqnarray*}
\hat\phi_{n,S} &=& \phi_{S}^*+
\Big(\frac{1}{n}\mathfrak{X}_{SS}\Big)^{-1}\left[\frac{1}{n}\sum_{t=1}^nZ_t{\bf X_t^S}-\lambda_n\lambda_{n,S} \, \text{sgn}(\phi_S^*)\right] \\
\lambda_n\hat\xi_{n,S^c} &=& \mathfrak{X}_{S^cS}\mathfrak{X}_{SS}^{-1}\left[\frac{1}{n}\sum_{t=1}^nZ_t{\bf X_t^S}-\lambda_n\lambda_{n,S} \, \text{sgn}(\phi_S^*)\right] -\frac{1}{n}\sum_{t=1}^n Z_t{\bf X_t^{S^c} } \; .
\end{eqnarray*}

Now, \emph{sign consistency} is equivalent (see Lemma 1 in \cite{wainwright})
to showing that
\begin{eqnarray}
\left|\phi_S^*+\big(\frac{1}{n}\mathfrak{X}_{SS}\big)^{-1}\Big[
\frac{1}{n}\sum_{t=1}^nZ_t{\bf X_t^S}-\lambda_n\lambda_{n,S}\, \text{sgn}
(\phi_S^*)\Big]\right|&>&0\label{eq:first} \\ \left| \mathfrak{X}_{S^cS}
\mathfrak{X}_{SS}^{-1}\Big[\frac{1}{n}\sum_{t=1}^nZ_t{\bf X_t^S}-\lambda_n
\lambda_{n,S}\, \text{sgn}(\phi_S^*)\Big] -\frac{1}{n}\sum_{t=1}^nZ_t{\bf X_t^{S^c}}
\right|&\leq& \lambda_n\lambda_{n,S^c}\Big] \label{eq:second}
\end{eqnarray}
hold, elementwise, with probability tending to $1$. Denote the
events in~(\ref{eq:first}), and in~(\ref{eq:second}) by
$\mathcal{A}$ and $\mathcal{B}$, respectively. The rest of the proof is devoted to showing that
$\mathbb{P}(\mathcal{A})\rightarrow1$, and
$\mathbb{P}(\mathcal{B})\rightarrow1$, as $n\rightarrow\infty$.

We commence with $\mathcal{A}$.  Let $\alpha_n=\min_{j\in S}|\phi_j^*|$.
Recall the notation $\|x_I\|_\infty$ for the $l_\infty$ norm on a set of
 indices $I$, i.e., $\max_{i\in I }|x_i|$ (and similarly for matrices).  It is enough to show that
 $\mathbb{P}(\|A_S\|_\infty>\alpha_n)\rightarrow 0$, as $n$ tends to infinity,  where
\begin{equation}\label{eq:A_S}
A_S = \Big(\frac{1}{n}\mathfrak{X}_{SS}\Big)^{-1}\Big[\frac{1}{n}
\sum_{t=1}^nZ_t{\bf X_t^S}-\lambda_n\lambda_{n,S}\, \text{sgn}(\phi_S^*)\Big] \; .
\end{equation}
Confine attention to the matrix $\mathfrak{X}_{SS}$.
The entry at row $i\in S$ and column $j\in S$ is given by
 $\sum_{t=1}^n X_{t-i} X_{t-j}$.  Notice that, equivalently,
  we can write this as $\sum_{t=1-i}^{n-i} X_t X_{t+i-j}$.
  Following~\cite{MR1093459}, one can show that $n^{-1}\mathfrak{X}_{SS}\rightarrow \Gamma_{SS}$
  in probability, as $n\rightarrow \infty$, where $\Gamma_{SS}=\big(\gamma(i-j)\big)_{i\in S, j\in S}$,
   and $\gamma(\cdot)$ is the autocovariance function, $\gamma(h)=\mathbb{E}X_tX_{t+h}$.
    Therefore, by assumption {\bf (i)} in Theorem \ref{thm:sign}, there exists a finite constant $C_{\max}$, such that $\|(n^{-1}\mathfrak{X}_{SS})^{-1}\|_\infty\leq o_P(1)+C_{\max}$.
     We continue by investigating the probability associated with the term
     inside the square brackets in~(\ref{eq:A_S}).

Notice that $\|\sum_{t=1}^n Z_t {\bf X_t^S }\|_\infty$ is given
by $\max_{i\in S}|\sum_{t=1}^n Z_tX_{t-i}|$, where $Z_t$ and
$X_{t-i}$ are independent random variables for each $t=1,
\ldots, n$, and $i\in S$. Fix an $i\in S$, and define 
\begin{equation}\label{eq:Tmartingale}
T_n\equiv
T_{n,i}=\sum_{t=1}^n Z_t X_{t-i} \; .
\end{equation}
Let
$\mathcal{F}_n=\sigma(\ldots, Z_{n-1}, Z_n)$ be the sigma-field
generated by $\{\ldots, Z_{n-1}, Z_n\}$. Simple calculation shows that
$\{T_n, \mathcal{F}_n\}_n$ is a martingale. Finally, Let
$Y_n=T_n-T_{n-1}$ denote the martingale difference sequence associated with
$T_n$. We quote below a  result concerning martingales
moment inequalities, which we shall make use of.
\begin{thm}[Burkholder's Inequality]
Let $\{X_n, \mathcal{F}_n\}_{n=1}^\infty$ be a martingale, and
$\tilde X_n=X_n-X_{n-1}$ be the associated martingale difference
sequence. Let $q>1$. For any finite and positive constants
$c=c(q)$, and $C=C(q)$ (depending only on $q$), we have
\begin{equation}\label{eq:burk}
c\Big[\mathbb{E}\big(\sum_{i=1}^n \tilde X_i^2\big)^{q/2}
\Big]^{1/q} \geq \big[\mathbb{E}|X_n|^q\big]^{1/q}\leq
C\Big[\mathbb{E}\big(\sum_{i=1}^n \tilde X_i^2\big)^{q/2}
\Big]^{1/q} \; .
\end{equation}
\end{thm}
Applying Cauchy-Schwartz inequality followed by Burkholder's
inequality, we obtain
\begin{equation}\label{eq:burk_applied}
\mathbb{E}|T_n| \leq \Big[\mathbb{E}\big|\sum_{t=1}^n Z_t
X_{t-i}\big|^2\Big]^{1/2} \leq C \Big[\sum_{t=1}^n
\mathbb{E}|Z_t^2 X_{t-i}^2|\Big]^{1/2}\leq C \sigma \sqrt{n} \; ,
\end{equation}
where $C$ is a finite and positive constant (from Burkholder's inequality). The last inequality
follows by the independence between $Z_t$ and $X_{t-i}$, and
since $\mathbb{E}|X_{t-i}|^2=1$. 
Fix an arbitrary, positive  $\xi<\infty$. By a trivial
bound we get
\begin{eqnarray*}
\mathbb{E}\max_{i\in S} |T_{n,i}|&\leq& \xi +\sum_{i\in
S}\int_\xi^\infty \mathbb{P}[|T_{n,i}|>y] \, dy \\
&\leq & \xi+\frac{1}{\xi}\sum_{i\in S}\mathbb{E}|T_{n,i}|^2 \\
&\leq& \xi+C^2\sigma^2 \frac{1}{\xi}sn \; ,
\end{eqnarray*}
recalling~(\ref{eq:burk_applied}). Now, picking $\xi=\sqrt{sn}$,
which is optimal, in the sense of obtaining an (asymptotically)
smallest fraction, we have,
\begin{equation}\label{eq:maximal_moment}
\frac{1}{n} \mathbb{E}\max_{i\in S} |T_{n,i}|\leq \sqrt{s/n}+ C^2\sigma^2\sqrt{s/n} = \OO{\sqrt{s/n}}
\; .
\end{equation}
This, in turn, implies, utilizing~(\ref{eq:A_S}) and Markov's
inequality, that $\mathbb{P}(\mathcal{A})\rightarrow 1$, by
imposing the condition:
\begin{equation*}
\frac{1}{\alpha_n} \Big[  \sqrt{s/n}  +
 \lambda_n\|\lambda_{n,S}\|_\infty\Big] \longrightarrow 0 \quad ,\quad \text{as} \; n\rightarrow\infty \; ,
\end{equation*}
which is condition (\ref{eq:sign_cond2}).

We turn to the event $\mathcal{B}$. Repeating the argument below~(\ref{eq:A_S}),
it is enough to show similar assertion about the event $\mathcal{B}$, with the
modification of replacing $\mathfrak{X}_{S^cS}\mathfrak{X}_{SS}^{-1}$,
by $\Gamma_{S^cS}\Gamma_{SS}^{-1}$. A sufficient condition for this to hold is
that $\{\| B_{S^c}\|_\infty\leq \lambda_n\min_{i\in S^c}\lambda_{n,i}\}$ happens with
probability tending to one, where
\begin{equation}\label{eq:B_Sc}
B_{S^c} = \Gamma_{S^cS}\Gamma_{SS}^{-1} \Big[\frac{1}{n}
\sum_{t=1}^nZ_t{\bf X_t^S}-\lambda_n\lambda_{n,S}\, \text{sgn}(\phi_S^*)\Big] -\frac{1}{n}\sum_{t=1}^nZ_t{\bf X_t^{S^c}}  \; .
\end{equation}

Under the \emph{incoherence condition} (condition {\bf (ii)} in the statement of the theorem), we have the following upper bound:
\[
\|B_{S^c}\|_\infty \leq (1-\epsilon)\frac{1}{n}\|
\sum_{t=1}^n Z_t{\bf X_t^S}\|_\infty+(1-\epsilon)\lambda_n \|\lambda_{n,S}\|_\infty +\frac{1}{n}\|
\sum_{t=1}^n Z_t{\bf X_t^{S^c}}\|_\infty \; ,
\]
which leads to: $\mathbb{P}(\|B_{S^c}\|_\infty>\lambda_n\min_{i\in
S^c} \lambda_{n,i}) \leq$
\begin{equation} \label{eq:prob_B}
\mathbb{P}\Big( \frac{2(1-\epsilon)}{n \lambda_n  \min_{i\in
S^c}\lambda_{n,i}}\| \sum_{t=1}^n Z_t{\bf X_t^S}\|_\infty> b \Big) +
\mathbb{P}\Big( \frac{2}{n \lambda_n
 \min_{i\in S^c}\lambda_{n,i}}\|\sum_{t=1}^n Z_t{\bf X_t^{S^c}}\|_\infty> b
 \Big)\; ,
\end{equation}
%\begin{eqnarray} \label{eq:prob_B}
%\lefteqn{\mathbb{P}(\|B_{S^c}\|_\infty>\lambda\min_{i\in S^c} \lambda_i) \leq }\\
%&& \nonumber \mathbb{P}\Big( \frac{2(1-\epsilon)}{n \lambda  \min_{i\in S^c}\lambda_i}\|
%\sum_{t=p+1}^n Z_t{\bf X_t^S}\|_\infty> b \Big) + \mathbb{P}\Big( \frac{2}{n \lambda
% \min_{i\in S^c}\lambda_i}\|\sum_{t=p+1}^n Z_t{\bf X_t^{S^c}}\|_\infty> b \Big) \; ,
%\end{eqnarray}
with $b=1-(1-\epsilon) \|\lambda_{n,S}\|_\infty / \min_{i\in S^c}\lambda_{n,i}$.
Note that inequality~(\ref{eq:prob_B}) follows by the inclusion
$\{U+V>z\}\subset \{U>z/2\}\cup \{V>z/2\}$. Under condition (\ref{eq:sign_cond1}), it would be enough to consider
the right hand side of~(\ref{eq:prob_B}), replacing (the two instances of) $b$  by
$\epsilon$. For the first term in~(\ref{eq:prob_B}) we have
\begin{equation}
\mathbb{P}\Big( \frac{2(1-\epsilon)}{n \lambda_n  \min_{i\in
S^c}\lambda_{n,i}}\| \sum_{t=1}^n Z_t{\bf X_t^S}\|_\infty> \epsilon \Big)\leq \frac{1-\epsilon}{\epsilon}\frac{2}{\lambda_n \min_{i\in S^c}\lambda_{n,i}} \frac{1}{n}\mathbb{E}\max_{i\in S^c} |T_{n,i}| \; ,
\end{equation}
which tends by~(\ref{eq:maximal_moment}) to zero once 
\begin{equation}\label{eq:cond3_1}
\frac{n\lambda_n ^2 (\min_{i\in S^c} \lambda_{n,i})^2}{s}\longrightarrow \infty \qquad , \qquad \text{as}\quad  n\rightarrow \infty \; .
\end{equation}
The same argument may be adapted for $\max_{i\in S^c}|T_{n,i}|$. We only need to replace $S$ by $S^c$. In this case we find that the condition
\begin{equation}\label{eq:cond3_2}
\frac{n\lambda_n ^2 (\min_{i\in S^c} \lambda_{n,i})^2}{\nu}\longrightarrow \infty \qquad , \qquad \text{as}\quad  n\rightarrow \infty \; ,
\end{equation}
is sufficient for showing that the second term in~(\ref{eq:prob_B}) converges to zero. 
Condition (\ref{eq:sign_cond3}) in the statement of the theorem guarantees both (\ref{eq:cond3_1}) and (\ref{eq:cond3_2}). The proof is now complete.

\end{proof}

\subsection{Estimation and Prediction Consistency}

\begin{proof}[Proof of Theorem \ref{thm:estimation}]

We follow~\cite{fan-2004-32}. In particular, denoting $\alpha_n=p^{1/2}(n^{-1/2}+\lambda_n\|\lambda_{n,S}\|)$, we will show that for every $\epsilon>0$ there exists a constant $C$, large enough, such that 
\begin{equation*}
\mathbb{P}\Big[\inf_{\|u\|=C}M_{\Lambda_n}(\phi^*+\alpha_n u)>M_{\Lambda_n}(\phi^*) \Big]>1-\epsilon \; ,
\end{equation*}
where $M_{\Lambda_n}(\cdot)$ is the objective function and is given in (\ref{eq:objective_function}).
This implies that $\|\hat\phi_n-\phi^*\|=O_P(\alpha_n)$. 

Multiplying both sides by $n$ clearly does not change the probability. We will show that $-n(M_{\Lambda_n}(\phi^*+\alpha_nu)-M_{\lambda_n}(\phi^*)) <0$ holds  uniformly over $\|u\|=C$. Write 
\[
M_{\Lambda_n}(\phi) = h(\phi) +\lambda_n \sum_{j=1}^p \lambda_{n,j} |\phi_j| \; ,
\]
for $h(\phi)=\|y-X\phi\|^2 / 2n$. We have
\begin{equation*}
-n(M_{\Lambda_n}(\phi^*+\alpha_nu)-M_{\Lambda_n}(\phi^*)) \leq -n[h(\phi^*+\alpha_nu)-h(\phi^*)]  -  n\lambda_n \sum_{j\in S} \lambda_{n,j} [|\phi_j^*+\alpha_nu_j|-|\phi_j^*|] \; .
\end{equation*}
Consider separately the least squares term, and the term associated with the $l_1$-penalty. We have, exploiting the fact that $\sum_{t=1}^n X_t {\bf X_t}=\mathfrak{X}\phi^*+\sum_{t=1}^n Z_t{\bf X_t}$, 
\begin{equation*}
-n[h(\phi^*+\alpha_nu)-h(\phi^*)]=\alpha_nu'\sum_{t=1}^n Z_t{\bf X_t}-\alpha_n^2u'\mathfrak{X}u/2 \equiv I_1-I_2 \; .
\end{equation*}
Recalling the definition of $T_{n,i}=\sum_{t=1}^n Z_tX_{t-i}$ (see (\ref{eq:Tmartingale})), and utilizing the result in~(\ref{eq:burk_applied}) we obtain
\begin{equation*}
|I_1|\leq \alpha_n\|u\| \| \sum_{t=1}^n Z_t{\bf X_t}\|=\|u\| O_P(\alpha_n \sqrt{pn}) \; .
\end{equation*}
Moving on to $I_2$, we write
\begin{equation}\label{eq:I2}
I_2=\alpha_n^2u'\mathfrak{X}u/2=n\alpha_n^2 u'(n^{-1}\mathfrak{X}-\Gamma_p)u/2 +
n \alpha_n^2u'\Gamma_p u/2 \; .
\end{equation}
We know that $n^{-1}\mathfrak{X_{ij}}$ tends in probability to $\gamma(i-j)$, where $\mathfrak{X}_{ij}=\sum_{t=1}^n X_{t-i}X_{t-j}$, the $(i,j)$ entry of $\mathfrak{X}$.  This clearly implies  $\|n^{-1}\mathfrak{X}-\Gamma_p\|=o_P(1)$, in the fixed $p$ scenario. Lemma~\ref{lem:consis_X} below shows that this may also hold true in the growing $p$ scenario which we consider here.
\begin{lem}\label{lem:consis_X}
Assume $\sum_{j=0}^\infty |\psi_j|<\infty$, as before. Then, 
\begin{equation}
\|n^{-1}\mathfrak{X}-\Gamma_p\|=o_P(1)  \; .
\end{equation}
\end{lem}
\begin{proof}
We adopt arguments given in~\cite[p. 226-227]{MR1093459}. Let $\epsilon>0$ be given. Using the fact that $\|A\|\leq \|A\|_F$, where $\|\cdot\|_F$ is the Frobenius matrix norm, $\sum_{i,j}|A_{ij}|^2$, we have
\begin{equation}\label{eq:consis_X}
\mathbb{P}(\|n^{-1}\mathfrak{X}-\Gamma_p\|>\epsilon)\leq \frac{1}{\epsilon^2}\sum_{i,j=1}^p d_{ij} \; ,
\end{equation}
where $d_{ij}=\mathbb{E}(n^{-1}\mathfrak{X}_{ij}-\gamma(i-j))^2$. We shall make use of Wick's formula. This formula gives the expectation of a product of several centered (joint) Gaussian variables $G_1, \ldots, G_N$, in terms of the elements of their covariance matrix $C=(c_{ij})$:
\begin{equation*}
\mathbb{E} \prod_{i=1}^k G_i=\sum c_{i_1i_2}\cdots c_{i_{k-1}i_k} \; ,
\end{equation*}
for $k=2m$, and zero otherwise. The sum extends over all different partitions of $\{G_1, \ldots, G_{2m}\}$ into $m$ pairs.
Applying the formula, we obtain:
\begin{eqnarray*}
\mathbb{E}\mathfrak{X}_{ij}^2&=&\sum_{s,t=1-i}^{n-i} 
\mathbb{E}X_tX_{t+i-j}X_sX_{s+i-j} \\ &=& \sum_{s,t=1-i}^{n-i} 
\Big(\gamma^2(i-j)+\gamma^2(s-t)+\gamma(s-t+i-j)\gamma(-(s-t)+i-j)\Big) \; ,
\end{eqnarray*}
where we have used the equivalent representation $\mathfrak{X}_{ij}=\sum_{t=1-i}^{n-i} X_t X_{t+i-j}$. 

A change of variables $k=s-t$ shows that
\begin{eqnarray*}
\lefteqn{\sum_{s,t=1-i}^{n-i} \Big(\gamma^2(s-t)+\gamma(s-t+i-j)\gamma(-(s-t)+i-j)\Big)  =} \\ && n[\gamma^2(0)+\gamma^2(i-j)] +2\sum_{k=1}^{n-1} (n-k)[\gamma^2(k)+\gamma(k+i-j)\gamma(-k+i-j)] \; .
\end{eqnarray*}
Therefore,
\begin{eqnarray}
d_{ij}=\lefteqn{\frac{p^2}{n^2}\gamma^2(i-j)+\frac{1}{n}[\gamma^2(0)+\gamma^2(i-j)] } \nonumber \\ \label{eq:d_ij} &+& \frac{2}{n^2} \sum_{k=1}^{n-1} (n-k)[\gamma^2(k)+\gamma(k+i-j)\gamma(-k+i-j)] \; .
\end{eqnarray}
Notice that $\sum_{k=1}^\infty |\gamma^2(k)+\gamma(k+i-j)\gamma(-k+i-j)|<\infty$. This may be seen by using the expression for the autocovariance function, $\gamma(h)=\sigma^2 \sum_{j=0}^\infty \psi_j\psi_{j+|h|}$, and by utilizing the summability of the $\psi_j$'s, $\sum_{j=0}^\infty |\psi_j| < \infty$. The expression~(\ref{eq:d_ij}) is therefore bounded by an $O(1/n)$ order term. This, in turn, shows that $d_{ij}=O(1/n)$, uniformly for every $i,j$. The proof is completed by recalling the RHS of~(\ref{eq:consis_X}), which is of the order of magnitude of $O(p^2/n)$.
\end{proof}

Using Lemma~\ref{lem:consis_X} we obtain 
\begin{equation}
|n\alpha_n^2 u'(n^{-1}\mathfrak{X}-\Gamma_p)u/2|\leq o_P(1)n\alpha_n^2 \|u\|^2 \; .
\end{equation}

We complete the argument with a bound on the term associated with the penalties, $-  n\lambda_n \sum_{j\in S} \lambda_{n,j} [|\phi_j^*+\alpha_nu_j|-|\phi_j^*|]$. Applying the Cauchy-Schwarz inequality, along with the fact that $\|a\|_1\leq \sqrt{p}\|a\|_2$ for every $a\in\mathbb{R}^p$, it is clear that the above term is absolutely bounded by $\lambda_n \|\lambda_{n,S}\|_\infty \sqrt{s}n\alpha_n \|u\|$. Now, since the second term in $I_2$ (see~(\ref{eq:I2})) dominates the other terms, the proof of the theorem is completed.

\end{proof}

\begin{proof}[Proof of Theorem \ref{thm:prediction}]

We begin as in \cite{buneaRandom}. Recall that $\|a\|_A^2$ stands for $a'Aa$, for every $p$-dimensional vector $a$, and $p\times p$ symmetric matrix $A$. We proceed by stating and proving two lemmas.
\begin{lem}\label{lem:aux1}
Let assumptions {\bf (i)}, and {\bf (ii)} of Theorem \ref{thm:prediction} be in effect. Then, 
\begin{equation}\label{eq:lemPred_1}
\|\hat\phi_n-\phi^*\|^2_{\mathfrak{X}/n} \leq 4\lambda_n M (s\kappa_p^{-1})^{1/2}\|\hat\phi_n-\phi^*\|_{\Gamma_p}
\end{equation}
holds true on
\begin{equation}\label{eq:event_1}
\mathcal{I}_1=\Big\{|\frac{2}{n}\sum_{t=1}^n X_{t-j}Z_t|\leq \lambda_n\lambda_{n,j} \quad , \quad \text{for all }\quad j=1,\ldots, p \Big\} \; .
\end{equation}

\end{lem}

\begin{proof}
By definition, the Lasso estimator $\hat \phi_n$ satisfies (see (\ref{eq:objective_function})),
\begin{equation*}
n^{-1} \|y-X\hat\phi_n\|^2+2 \lambda_n \sum_{j=1}^p\lambda_{n,j}|\hat\phi_{n,j}| \leq 
n^{-1} \|y-X\phi^*\|^2+2\lambda_n \sum_{j=1}^p\lambda_{n,j}|\phi^*_j| \; .
\end{equation*}
Recalling the model $y=X\phi^*+Z$, we obtain, by re-arrangeing the above terms,
\begin{equation*}
\|\hat\phi_n-\phi^*\|^2_{\mathfrak{X}/n}+2\lambda_n\sum_{j=1}^p\lambda_{n,j}|\hat\phi_{n,j}| \leq 2(\hat\phi_n-\phi^*)'\frac{1}{n}X'Z +2\lambda_n\sum_{j=1}^p\lambda_{n,j}|\phi^*_j| \; .
\end{equation*}
Now, since $(\hat\phi_n-\phi^*)'\frac{1}{n}X'Z=\sum_{j=1}^p (\hat\phi_{n,j}-\phi^*_j)\frac{1}{n}\sum_{t=1}^n X_{t-j}Z_t$, we have, on $\mathcal{I}_1$,
\begin{eqnarray}
\|\hat\phi_n-\phi^*\|^2_{\mathfrak{X}/n}&\leq& \lambda_n \sum_{j=1}^p\lambda_{n,j} |\hat\phi_{n,j}-\phi^*_j| +2\lambda_n \sum_{j=1}^p \lambda_{n,j} (|\phi_j^*|-|\hat\phi_{n,j}|) \nonumber   \\ 
&\leq & 4\lambda_n \sum_{j\in S}\lambda_{n,j} |\hat\phi_{n,j}-\phi^*_j| \; , \label{eq:lem_diff}
\end{eqnarray}
where the second inequality is obtained by decomposing  the summation $\sum_{j=1}^p$ into $\sum_{j\in S}+\sum_{j\notin S}$, and using Cauchy-Schwarz inequality.

By assumption {\bf (ii)}, and the fact that $\gamma(0)=\mathbb{E}|X_t|^2=1$, we have
\begin{eqnarray}
\nonumber 
\sum_{j\in S}|\hat\phi_{n,j}-\phi^*_j|^2 &\leq& \sum_{j=1}^p(\hat\phi_{n,j}-\phi_j^*)^2 = \|\hat\phi_n-\phi^*\|^2_{\text{diag}(\Gamma_p)} \nonumber \\
&\leq& \frac{1}{\kappa_p}\|\hat\phi_n-\phi^*\|^2_{\Gamma_p} \; .  \label{eq:diagGamma}
\end{eqnarray}
The proof is completed by applying the Cauchy-Schwarz inequality on (\ref{eq:lem_diff}), and by using assumption {\bf (i)}.

\end{proof}

We turn to the second lemma.
\begin{lem}\label{lem:aux2}
Let assumptions {\bf (i)}, {\bf (ii)} of Theorem \ref{thm:prediction} be in effect. Let $C$ be a constant (given explicitly in the proof) depending on $M$ only. Put $\epsilon = \lambda_n (sp^{-1})^{1/2}$. Then,
\begin{equation}
\|\hat\phi_n-\phi^*\|^2_{\Gamma_p} \leq C \lambda_n^2 s \kappa_p^{-1} \; ,
\end{equation}
holds true on $\mathcal{I}_1\cap \mathcal{I}_2$, where $\mathcal{I}_1$ is given by (\ref{eq:event_1}), and 
\begin{equation}
\mathcal{I}_2=\left\{M_p \leq \epsilon \right\} \; ,
\end{equation}
with
\begin{equation}
M_p  = \max_{1\leq i,j \leq p} \left|\frac{\mathfrak{X}_{ij}}{n}-\gamma(i-j)\right| \; .
\end{equation}
\end{lem}

\begin{proof}
Note that 
\[
\left| \|\hat\phi_n-\phi^* \|^2_{\mathfrak{X}/n} - \| \hat\phi_n-\phi^*\|^2_{\Gamma_p} \right| \leq M_p \|\hat\phi_n-\phi^*\|^2_1 \; .
\]
Therefore,
\begin{eqnarray*}
\|\hat\phi_n-\phi^*\|^2_{\mathfrak{X}/n}&\geq& \|\hat\phi_n-\phi^*\|^2_{\Gamma_p}-M_p p^{1/2}\|\hat\phi_n-\phi^*\|  \\
&\geq& \|\hat\phi_n-\phi^*\|^2_{\Gamma_p}-M_p (p\kappa_p^{-1})^{1/2}\|\hat\phi_n-\phi^*\|_{\Gamma_p} \; .
\end{eqnarray*}
The first inequality follows since $\|a\|_1\leq n\|a\|^2$, and the second inequality is satisfied under assumption {\bf (ii)} (see (\ref{eq:diagGamma})).
Referring back to (\ref{eq:lemPred_1}), we obtain, on $\mathcal{I}_1\cap \mathcal{I}_2$,
\[
\|\hat\phi_n-\phi^*\|^2_{\Gamma_p} \leq 2 (1/2+2M)\lambda_n (s\kappa_p^{-1})^{1/2} \|\hat\phi_n-\phi^*\|_{\Gamma_p}  \; .
\]
Applying the inequality $2xy\leq 2x^2+y^2/2$ on the right-hand side of the expression above (with $x=(1/2+2M) \lambda_n (s\kappa_p^{-1})^{1/2}$, and $y=\|\hat\phi_n-\phi^*\|_{\Gamma_p}$), we establish the statement of the Lemma, with $C=4(1/2+2M)^2$.

\end{proof}

The rest of the proof of Theorem \ref{thm:prediction} is devoted to showing that indeed  $\|\hat\phi_n-\phi^*\|^2_{\Gamma_p}\leq C\lambda_n^2 sk_p^{-1}$ holds on a negligible event, i.e., that the probability of the complement of $\mathcal{I}_1\cap \mathcal{I}_2$ is negligible. We shall commence with $\mathcal{I}_2$.

We recall here the family of time series $\{X_t\}$, denoted by $\mathcal{H}_\rho(l,L)$, for some $\rho>1$, $0<l<1$, and $L>1$ (Section \ref{sec:estimation_prediction}). The family consists of all stationary Gaussian time series with $\mathbb{E}X_t=0$, $\mathbb{E}|X_t|^2=1$, and 
enjoys an exponential decay of the strongly mixing coefficients (see (\ref{eq:decay_alpha})).

\begin{lem}\label{lem:I2}
Assume that $\epsilon=\lambda_n(s/p)^{1/2}\leq Dn^{-2/5}$, where  $D=(C_1^3C_2\beta_1^2\beta_2^3)^{1/5}$, with $C_1$ and $C_2$ two constants explicitly specified in the proof. Then,
\[
\mathbb{P}(\mathcal{I}_2^c)\leq p^2 \exp\big\{-n\lambda_n^2 (s/p^2) / (4C_1 \beta_1\beta_2)\big\} \; .
\]
\end{lem}

\begin{proof}
We begin with 
\[
\mathbb{P}(|\sum_{t=1-i}^{n-i} Y_t|>\epsilon) \; ,
\]
where
\begin{equation}\label{eq:Y}
Y_t \equiv Y_{t,i,j} = \frac{1}{n}(X_tX_{t+i-j}-\gamma(i-j)) \; .
\end{equation}
The proof is based on an application of the pair of lemmas \ref{lem:zeevi1} and \ref{lem:zeevi2}, after noticing that
\[
\mathbb{	P}(\mathcal{I}_2^c)=\mathbb{P}(M_p>\epsilon)\leq \sum_{i,j=1}^p \mathbb{P}\left(\big|\sum_{t=1-i}^{n-i} Y_t\big|>\epsilon\right) \; .
\]
Define $k=i-j$. It is enough to consider only $k\geq 0$ ($i\geq j$), since $\mathfrak{X}_{ij}$ and $\gamma(i-j)$ are symmetric. By the same argument below expression (39) in \cite{zeevi}, one may notice that $\{Y_t\}$ is strongly mixing with the rate $\alpha_Y(m)\leq \alpha_X(m-k)$ for $m>k$, and $\alpha_Y(m)\leq 1/4$ (see \cite{bradley}), but for our purposes in would be enough to bound $\alpha_Y(m)$, for $m>k$, by simply $1$.

We shall make use of the following two lemmas, adapted from \cite{zeevi}.

\begin{lem}\label{lem:zeevi1}
Suppose $\{X_t\}$ is a strongly mixing time series, $S_n=\sum_{t=1}^n$, and $\text{cum}_r(S_n)$ is the $r$th order cumulant of $S_n$. For $\nu>0$ define the function
\[
\Lambda_n(\alpha_X, \nu) = \max \big\{ 1 \, , \, \sum_{m=1}^n (\alpha_X(m))^{1/\nu} \big\} \; .
\]
If, for some $\mu\geq 0$, $H>0$
\[
\mathbb{E}|X_t|^r \leq (r!)^{\mu+1} H^r \quad , \quad t=1,\ldots , n  , \, r=2,3,\ldots  \; ,
\]
then $|cum_r(S_n)|\leq 2^{r(1+\mu)+1}12^{r-1}(r!)^{2+\mu} H^r [\Lambda_n(\alpha_X, 2(r-1))]^{r-1} n$.
\end{lem}

\begin{lem}\label{lem:zeevi2}
Let $Y$ be a random variable with $\mathbb{E} Y=0$. If there exist $\mu_1\geq 0$, $H_1>0$ and $\Delta>0$ such that
\[
|cum_r(Y)|\leq \left(\frac{r!}{2}\right)^{1+\mu_1} \frac{H_1}{\Delta^{r-2}} \; , \quad  r=2,3,\ldots \; ,
\]
then
\[
\mathbb{P}(|Y|>y) \leq \left\{
\begin{array}{ll}
\exp\{-y^2/(4H_1)\} & 0\leq y\leq (H_1^{1+\mu_1}\Delta)^{1/(2\mu_1+1)} \\
\exp\{-(y\Delta)^{1/(1+\mu_1)}/4\} & y\geq (H_1^{1+\mu_1}\Delta)^{1/(2\mu_1+1)} \; .
\end{array}
\right.
\]
\end{lem}

Back to the proof of Lemma \ref{lem:I2}. Absolute moment of $Y_t$ are bounded as follows:
\begin{eqnarray*}
\mathbb{E}|Y_t|^r &\leq& n^{-r} 2^{r-1} \big[\mathbb{E}|X_tX_{t+k}|^r+|\gamma(k)|^r\big] \\
&\leq & n^{-r}2^{r-1} \left[ \big(\mathbb{E}|X_t|^{2r}\mathbb{E}|X_{t+k}|^{2r}\big)^{1/2}+\gamma(0) \right] \\
&\leq & r! (4/n)^r \; .
\end{eqnarray*}
The second inequality follows by the Cauchy-Schwarz inequality together with the inequality $(a+b)^j\leq 2^{j-1}(a^j+b^j)$, and the last inequality follows by the assumed Gaussianity of $X_t$, and the inequality ${2r \choose r}\leq 2^{2r}$.
We have
\begin{eqnarray*}
\sum_{m=1}^n (\alpha_X(m))^{1/2(r-1)} &\leq & k + \left(\frac{2Ll}{l(\rho-1)}\right)^{1/(r-1)} \sum_{m=1}^{n-k} \rho^{-m/2(r-1)} \\
&\leq &  k+ \left(\frac{2Ll}{l(\rho-1)}\right)^{1/(r-1)} \left(1+\frac{2(r-1)}{\log \rho} \right) \; ,
\end{eqnarray*}
The first inequality utilizes the relationship between $\alpha_Y(m)$ and $\alpha_X(m)$, and inequality (\ref{eq:decay_alpha}). The second inequality uses geometric series expression together with the inequality $\rho^x-1\geq x\log \rho$, for all $x\geq 0$.

Therefore, defining $\hat k=k$ if $k>0$, and $\hat k=1$, if $k=0$, we obtain, after some manipulations, similar to those in \cite{zeevi},
\[
[\Lambda_n(\alpha_X, 2(r-1))]^{r-1}\leq 12^{r-1} r! (\hat k \beta_1)^{r-1} \beta_2 \; ,
\]
for two constants $\beta_1$ and $\beta_2$, given, respectively, by $1+1/\log\rho$ and $1+L\rho/l(\rho-1)$ (see (\ref{eq:constants})). The bound results from the inequalities $(a+b)^j\leq 2^{j-1}(a^j+b^j)$, $n^n\leq n!e^{n}$, and other trivial inequalities.

Applying Lemma \ref{lem:zeevi1} (with $\mu=0$ and $H=4/n$) we have $|cum_r(\sum_{t=1-i}^{n-i})Y_t|\leq \text{RHS}$, where RHS can be put in the form $(r!/2)^3 H\Delta^{2-r}$, with $H_1=C_1\beta_2(\hat k\beta_1/n)$, $\Delta=C_2(\hat k \beta_1/n)^{-1}$, and $C_1=2^{10}12^2$, $C_2=2^{-3}12^{-2}$. Now, applying Lemma \ref{lem:zeevi2} (with $\mu_1=2$, and $H_1$ and $\Delta$ as above) we obtain:
\begin{equation}
\mathbb{P}\left(\big|\sum_{t=1-i}^{n-i} Y_t\big|>y\right) \leq \left\{
\begin{array}{ll} \label{eq:Y_large_dev}
\exp\left\{-y^2n/(4C_1\hat k \beta_1\beta_2)\right \} & 0\leq y\leq D\hat k^{2/5} n^{-2/5} \\
\exp\left\{-\frac{1}{4}\big(\frac{C_2}{\hat k \beta_1}\big)^{1/4}(yn)^{1/3}\right\} & y\geq D \hat k^{2/5} n^{-2/5} \; ,
\end{array}
\right.
\end{equation}
where $D=(C_1^3C_2 \beta_1^2\beta_2^3)^{1/5}$. The proof is completed by applying the moderate deviation part in (\ref{eq:Y_large_dev}) with $y=\epsilon$, and by noticing that $1\leq \hat k \leq p$.

\end{proof}

We turn to evaluate the probability of the complement of the event $\mathcal{I}_1$.

\begin{lem}\label{lem:I1}
For all $0<c<\infty$ and $y>\sigma^2(n+Dn^{3/5})$ (where $D$ is given by (\ref{eq:constants})),
\[
\mathbb{P}(\mathcal{I}_1^c) \leq 6p\exp\left\{ -F_1\min \left\{(\sigma^{-2}y-n)^{1/3}, c^2\sigma^{-2}, \frac{n^2\lambda_n^2\lambda^2_{\min}}{y+cn\lambda_n\lambda_{\max}/2}\right\}  \right\}  \; ,
\]
where $F_1=\min\left\{ (C_2 / \beta_1)^{1/4}/4, 2^{-9}, 8^{-1} \right\}$.
\end{lem}

\begin{proof}
Let $V_n^2 = \sigma^2 \sum_{t=1}^n X_{t-i}^2=\sigma^2\sum_{t=1-i}^{n-i} X_t^2$. 
Fix a $y>\sigma^2(n+Dn^{3/5})$ and a $0 < c < \infty$. Denote by $\tilde{\mathcal{I}}_1$ the event $\mathcal{I}_1$ (see (\ref{eq:event_1})) with the absolute value removed.  We begin by writing:
\begin{eqnarray*}
\mathbb{P}(\tilde{\mathcal{I}}_1^c) &\leq& \sum_{j=1}^p \mathbb{P}\big(\frac{2}{n}\sum_{t=1}^n X_{t-j}Z_t > \lambda_n\lambda_{n,j}\big) \\
&\leq & \sum_{j=1}^p \mathbb{P}\Big(\bigcup_{n=1}^\infty \big\{\frac{2}{n}\sum_{t=1}^n X_{t-j}Z_t >\lambda_n\lambda_{n,j} \, , \, V_n^2\leq y \big\} \Big) +p\, \mathbb{P}(V_n^2>y) \\
&=:& I_1 + I_2  \; .
\end{eqnarray*}
Clearly, $I_1$ satisfies $I_1 \leq I_{11}+I_{12}$, with
\begin{eqnarray*}
I_{11} &=&  \sum_{j=1}^p \mathbb{P}\Big(\bigcup_{n=1}^\infty \big\{\frac{2}{n}\sum_{t=1}^n X_{t-j}Z_t >   \lambda_n\lambda_{n,j} \, , \, V_n^2\leq y \big\} \, , \, \bigcap_{r=3}^\infty \big\{|X_{t-j}|^{r-2}\mathbb{E}|Z_t|^r\leq \frac{r!}{2}\sigma^2c^{r-2} \big\}\Big)  \; , \\
I_{12} &=& \sum_{j=1}^p \mathbb{P}\Big( \bigcup_{r=3}^\infty \big\{ |X_{t-j}|^{r-2}\mathbb{E}|Z_t|^r> \frac{r!}{2}\sigma^2c^{r-2}  \big\} \Big) \; .
\end{eqnarray*}
We analyze $\mathbb{P}(\tilde{\mathcal{I}}_1^c)$ by investigating $I_{11}$, $I_{12}$ and $I_2$ separately. 

For $I_2$, we recall that $Y_t\equiv Y_{t,i,i}=(X_{t}^2-\gamma(0))/n$ (see (\ref{eq:Y}) and the remark  below) is strongly mixing with exponential decay rate. Therefore, by the large deviation part in (\ref{eq:Y_large_dev}) (with $\hat k=1$),
\begin{eqnarray*}
\mathbb{P}(V_n^2>y) &\leq& \mathbb{P}(|V_n^2-n\sigma^2|>y-n\sigma^2) \\
&=& \mathbb{P}(|\sum_{t=1-i}^{n-i} Y_t|> \sigma^{-2}n^{-1} y -1) \\
&\leq & \exp\left\{-\frac{1}{4}\big(\frac{C_2}{\beta_1}\big)^{1/4}(\sigma^{-2}y-n)^{1/3}\right\} \; .
\end{eqnarray*}

For $I_{12}$, we use the bound $\mathbb{E}|Z_t|^{2r} \leq \sigma^{2r} r!2^{2r}$ (and the Cauchy-Schwarz inequality) to obtain 
\[
\big\{ |X_{t-j}|^{r-2}\mathbb{E}|Z_t|^r> \frac{r!}{2}\sigma^2c^{r-2}  \big\} \subset \big\{ |X_{t-j}| >2^{-(1+r)/(r-2)} \sigma^{-1}c \big\} \; .
\]
Therefore, noticing that $\{2^{-(1+r)/(r-2)}\}_{r=3}^\infty$ is an increasing sequence, we have
\[
I_{12} \le \sum_{j=1}^p \mathbb{P}\big( |X_{t-j}|> 2^{-4}\sigma^{-1}c\big ) \leq (2/\pi)^{1/2}  p \exp\{-2^{-8} c^2 / 2\sigma^2 \}  \; .
\]

For $I_{11}$, we use the following theorem which is a Bernstein's type of an inequality for martingales.  
\begin{thm}[\citet{pena}]
Let $\{M_n, \mathcal{F}_n\}$ be a martingale, with difference $\Delta_n = M_n-M_{n-1}$. Define $V_n^2=\sum_{i=1}^n \sigma_i^2=\sum_{i=1}^n \mathbb{E}(\Delta_i^2 \, | \, \mathcal{F}_{i-1})$. Assume that $\mathbb{E}(|\Delta_i|^r \, | \, \mathcal{F}_{i-1})\leq (r!/2)\sigma_i^2c^{r-2}$ a.e. for $r\geq 3$, $0<c<\infty$. Then, for all $x,y>0$,
\begin{equation}\label{eq:pena}
\mathbb{P}\Big(\bigcup_{n=1}^\infty \{ M_n> x \, , \, V_n^2 \leq y \} \Big) \leq \exp\left\{ -\frac{x^2}{2(y+cx)} \right\} \; .
\end{equation}
\end{thm}
\noindent Recall that $\sum_{t=1}^n X_{t-j}Z_t$ is a martingale (see (\ref{eq:Tmartingale})). Then, simple application of the above theorem, with $x=n\lambda_n\lambda_{n,j}/2$, leads to
\[
I_{11} \leq  p \exp\left\{  - \frac{n^2\lambda_n^2\lambda_{\text{min}}^2}{8(y+cn\lambda_n\lambda_{\text{max}}/2)} \right\} \; .
\]

Lemma \ref{lem:I1} now follows by collecting the bounds of $I_{11}$, $I_{12}$, and $I_2$, and by symmetry.

\end{proof}

The proof of theorem \ref{thm:prediction} is now complete by virtue of Lemma \ref{lem:aux1}, Lemma \ref{lem:aux2}, Lemma \ref{lem:I2}, and Lemma \ref{lem:I1}.

\end{proof}

%\bibliographystyle{plainnat}
%\bibliography{nlts}

\end{document}